\documentclass{amsart}

\usepackage{amssymb,latexsym,amsfonts,amsmath}
\usepackage{graphicx,psfrag,epsfig,epsf}
\include{diagrams}
\usepackage{eso-pic}
\usepackage{graphicx}
\usepackage{color}
\usepackage{diagrams}
\usepackage{algorithm2e}
\usepackage{tikz}
\usetikzlibrary{automata}

\usepackage{amsfonts}
\usepackage{amssymb,latexsym,amsfonts,amsmath}
\usepackage{graphicx}

\newtheorem{theorem}{Theorem}[section]
\newtheorem{lemma}[theorem]{Lemma}

\newtheorem{problem}[theorem]{Problem}
\newtheorem{proposition}[theorem]{Proposition}

\theoremstyle{definition}
\newtheorem{definition}[theorem]{Definition}

\theoremstyle{remark}

\numberwithin{equation}{section}

\begin{document}

\begin{abstract}                          % Abstract of not more than 200 words.
Discrete abstractions of continuous and hybrid systems have recently been the topic of great interest from both the control systems and the computer science communities, because they provide a sound mathematical framework for analysing and controlling embedded systems.
In this paper we give a further contribution to this research line, by addressing the problem of symbolic control design 
of nonlinear systems with infinite states specifications, modelled by differential equations. We first derive the symbolic controller solving the control design problem, given in terms of discrete abstractions of the plant and the specification systems. We then present an algorithm which integrates the construction of the discrete abstractions with the design of the symbolic controller. Space and time complexity analysis of the proposed algorithm is performed and a comparison with traditional approaches currently available in the literature for symbolic control design, is discussed. Some examples are included, which show the interest and applicability of our results.
\end{abstract}

\title[Integrated symbolic control design for nonlinear systems with infinite states specifications]{Integrated symbolic control design for nonlinear systems with infinite states specifications}
\thanks{This work has been partially supported by the Center of Excellence for Research DEWS, University of L'Aquila, Italy.}

\author[Giordano Pola, Alessandro Borri and Maria D. Di Benedetto]{
Giordano Pola$^{1}$, Alessandro Borri$^{1}$ and Maria D. Di Benedetto$^{1}$}
\address{$^{1}$
Department of Electrical and Information Engineering, Center of Excellence DEWS,
University of L{'}Aquila, Poggio di Roio, 67040 L{'}Aquila, Italy}
\email{ \{giordano.pola,alessandro.borri,mariadomenica.dibenedetto\}@univaq.it}
%\urladdr{
%http://www.diel.univaq.it/people/pola/
%}
%\urladdr{
%http://www.diel.univaq.it/people/pepe/
%}
%\urladdr{
%http://www.diel.univaq.it/people/dibenedetto/
%}

%\address{$^{2}$Department of Electrical Engineering\\
%University of California at Los Angeles,
%Los Angeles, CA 90095}
%\email{tabuada@ee.ucla.edu}
%\urladdr{http://www.ee.ucla.edu/$\sim$tabuada}

\maketitle

\section{Introduction}\label{sec1}

Discrete abstractions of continuous and hybrid systems have been the topic of intensive study in the last twenty years from both the control systems and the computer science communities \cite{TACsymbolicmodels}. While physical world processes are often described by differential equations, digital controllers and software and hardware at the implementation layer, are usually modelled through discrete/symbolic processes. This mathematical models heterogeneity has posed during the years interesting and challenging theoretical problems that are needed to be addressed, in order to ensure the formal correctness of control algorithms. One approach to deal with this heterogeneity is to construct symbolic models that are equivalent to the continuous process, so that the mathematical model of the process, of the controller, and of the software and hardware at the implementation layer, are of the same nature.
Several classes of dynamical and control systems admitting symbolic models, were identified during the years. We recall timed automata \cite{TheoryTA}, rectangular hybrid automata \cite{WhatsDecidable}, and o-minimal hybrid systems \cite{OMinimal} in the class of hybrid automata. Control systems were considered further. Early results in this regard are reported in the work of \cite{caines}, \cite{moor}, \cite{forstner} and \cite{BMP02}. %Recent results on symbolic models of control systems i
Recent results %on symbolic models of control systems
include the work of \cite{LTLControl}, which showed existence of symbolic models for controllable discrete--time linear systems, and the work of \cite{HCS06,Belta:06} for piecewise--affine and multi--affine systems. Many of the aforementioned work are based on the notion of bisimulation equivalence, introduced by Milner and Park \cite{Milner,Park} in the context of concurrent processes, as a formal equivalence notion to relate continuous and hybrid processes to purely discrete/symbolic models. A new insight in the construction of symbolic models has been recently placed through the notion of approximate bisimulation introduced by Girard and Pappas in \cite{AB-TAC07}. Based on the above notion, some classes of incrementally stable \cite{IncrementalS} control systems were recently shown to admit symbolic models: discrete--time linear control systems \cite{Girard_HSCC07}, nonlinear control systems with and without disturbances \cite{PolaAutom2008,PolaSIAM2009}, nonlinear time--delay systems \cite{PolaSCL10} and switched nonlinear systems \cite{GirardTAC2010}. Recent results in the work of \cite{MajidACC10} have also shown the existence of symbolic models for unstable nonlinear control systems, satisfying the so--called incremental forward completeness property.\\
The use of symbolic models in the control design of continuous and hybrid systems has been investigated in the work of \cite{LTLControl,BeltaCDC09,TabuadaTAC08}, among many others. The work in \cite{LTLControl} considers discrete--time linear control systems, the work in \cite{BeltaCDC09} considers piecewise--affine systems while the work in \cite{TabuadaTAC08} considers stabilizable nonlinear control systems. %The main drawback of using simbolic models for control design of continuous processes is the size of the symbolic models which often grows exponentially with the dimension of the state space of the continuous process to--be--controlled.
In this paper we give a further contribution to this research line and in particular, in the direction of \cite{TabuadaTAC08}.
We consider symbolic control design of nonlinear control systems where specifications are characterized by an infinite number of states and modelled through differential equations: Given a plant nonlinear control system and a specification nonlinear (autonomous) system, we investigate conditions for the existence of a symbolic controller that implements the behaviour of the specification, with a precision that can be rendered as small as desired. In other words, we look for a symbolic controller so that the interconnection between the plant and the controller satisfies or conforms \cite{ModelChecking} the specification with an arbitrarily small precision. The symbolic controller is furthermore requested to be non--blocking in order to prevent the occurrence of deadlocks in the interaction between the plant and the symbolic controller. %To the best of the authors knowledge, symbolic control design of systems with infinite states specifications has not been addressed before in the current literature.
This control design problem can be seen as an approximated version of similarity games, as discussed in \cite{paulo}. Similar problems have been studied in the literature (in a non--approximating settings) in the context of supervisory control \cite{CassandrasBook}, symbolic control design for piecewise--affine systems enforcing temporal logic specifications \cite{BeltaCDC09}, among many others.\\
%, see also classical control problems in the context of supervisory control, see e.g. \cite{CassandrasBook}.
The control design problem that we consider in this paper has been solved by following the so--called correct--by--design approach, see e.g. \cite{LTLControl,TabuadaTAC08,BeltaCDC09}. We first construct the symbolic models of the plant and the specification by making use of (some variations of) the results established in \cite{PolaAutom2008}. We then solve the control design problem at the symbolic layer, to finally come back at the continuous layer, by providing appropriate approximating bounds in the quantization errors which guarantee the solution to the control design problem under study. The solution of the control design problem at the symbolic layer is shown to be the maximal non--blocking part of the (exact) parallel composition \cite{CassandrasBook} of the symbolic models associated with the plant and the specification. By following the correct--by--design approach, the design of the symbolic controller solving the problem at hand, requires a first computation of the plant and the specification symbolic models, then a construction of the (exact) parallel composition of the symbolic systems obtained and finally a computation of the maximal non--blocking part of the composed system. 
While being formally correct from the theoretical point of view, this approach is in general rather demanding from the computational point of view, because of the large size of the symbolic models needed to be constructed, in order to synthesize the symbolic controller solving the design problem. This drawback is common with other approaches currently available in the literature on symbolic control design of continuous and hybrid systems, see e.g. \cite{LTLControl,BeltaCDC09,TabuadaTAC08} and motivated some researchers to propose solutions to cope with complexity. For example, the work in \cite{Tazaki09} proposes nonuniform state quantizations in the construction of the symbolic models of the to--be--controlled plant system. In this paper we propose an alternative solution to the one studied in \cite{Tazaki09}. Inspired by on--the--fly verification and control of timed or untimed transition systems (see e.g. \cite{onthefly3,onthefly2}), we approach the design of symbolic controllers by advocating an ``integration'' philosophy: instead of computing separately the symbolic models of the plant and of the specification to then design the controller at the symbolic layer, \textit{we integrate each step of the procedure in only one algorithm}.
Space and time complexity analysis of the proposed algorithm is performed and a comparison with traditional approaches currently available in the literature, is discussed. Some examples are included which show the interest and applicability of our results. 
%The proposed algorithm is proved by a space and time complexity analysis to perform better than traditional approaches currently available in the literature. Some examples are included which show the interest in and the applicability of the proposed results. 
For the sake of completeness, a detailed list of the employed notation is included in the Appendix (Section \ref{appendix}).

\section{Preliminary Definitions}\label{sec2}

\subsection{Control Systems\label{II.B}}

The class of control systems that we consider in this paper is formalized in
the following definition.

\begin{definition}
\label{Def_control_sys}A \textit{control system} is a quintuple:
\begin{equation}
\Sigma=(X,X_{0},U,\mathcal{U},f),
\label{NLCsystem}
\end{equation}
where:

\begin{itemize}
\item $X\subseteq\mathbb{R}^{n}$ is the state space;

\item $X_{0}\subseteq X$ is the set of initial states;

\item $U\subseteq\mathbb{R}^{m}$ is the input space;

\item $\mathcal{U}$ is a subset of the set of all locally essentially bounded
functions of time from intervals of the form \mbox{$]a,b[\subseteq\mathbb{R}$}
to $U$ with $a<0$, $b>0$;

\item \mbox{$f:\mathbb{R}^{n}\times U \rightarrow\mathbb{R}^{n}$} is a
continuous map satisfying the following Lipschitz assumption: for every
compact set \mbox{$K\subset\mathbb{R}^{n}$}, there exists a constant $\kappa\in\mathbb{R}^{+}$
such that 
\[
\Vert f(x,u)-f(y,u)\Vert\leq \kappa\Vert x-y\Vert,
\]
for all $x,y\in K$ and all $u\in U$.
\end{itemize}
\end{definition}

A curve \mbox{$\xi:]a,b[\rightarrow\mathbb{R}^{n}$} is said to be a
\textit{trajectory} of $\Sigma$ if there exists $u\in\mathcal{U}$
satisfying:
\begin{equation}
\dot{\xi}(t)=f(\xi(t),u(t)),
\label{eq0}%
\end{equation}
for almost all $t\in$ $]a,b[$. Although we have defined trajectories over open
domains, we shall refer to trajectories
\mbox{${\xi:}[0,\tau]\rightarrow\mathbb{R}^{n}$} defined on closed domains
$[0,\tau],$ $\tau\in\mathbb{R}^{+}$ with the understanding of the existence of
a trajectory \mbox{${\xi}^{\prime}:]a,b[\rightarrow\mathbb{R}^{n}$} such that
\mbox{${\xi}={\xi}^{\prime}|_{[0,\tau]}$}. We also write $\xi_{xu}%
(\tau)$ to denote the point reached at time $\tau$ under the input $u$
from initial condition $x$; this point is uniquely determined, since the
assumptions on $f$ ensure existence and uniqueness of trajectories
\cite{Sontag}. A control system $\Sigma$ is said to be forward complete if every trajectory is defined on an interval of the form $]a,\infty\lbrack$. Sufficient and necessary conditions for a system to be forward complete can be found in \cite{Angeli&Sontag1999}. The above formulation of control systems can be also used to model autonomous nonlinear systems, i.e. systems with no control inputs. With a slight abuse of notation we denote an autonomous system $\Sigma$ by means of the tuple $(X,X_{0},f)$. 

\subsection{Systems}

We will use systems to describe both control systems as well as their symbolic models. For a detailed exposition of the notion of systems and of their properties we refer to \cite{paulo}.

\begin{definition}
\cite{paulo} A system $S$ is a sextuple:
\[
S=(X,X_{0},U,\rTo,Y,H),
\]
consisting of:

\begin{itemize}
\item a set of states $X$;

\item a set of initial states $X_{0}\subseteq X$;

\item a set of inputs $U$;

\item a transition relation $\rTo\subseteq X\times U\times X$;

\item an output set $Y$;

\item an output function $H:X\rightarrow Y$.

\end{itemize}
\end{definition}

A transition $(x,u,x')\in\rTo$ of system $S$ is denoted by $x\rTo^{u}x^{\prime}$. System $S$ is said to be:

\begin{itemize}
\item \textit{countable}, if $X$ and $U$ are countable sets;

\item \textit{symbolic}, if $X$ and $U$ are finite sets;

\item \textit{metric}, if the output set $Y$ is equipped with a metric
$d:Y\times Y\rightarrow\mathbb{R}_{0}^{+}$;

\item \textit{deterministic}, if for any $x\in X$ and $u\in U$ there exists at most one $x^{\prime}\in X$ such that $(x,u,x^{\prime})\in\rTo$;

\item \textit{non--blocking}, if for any $x\in X$ there exists $(x,u,x^{\prime})\in\rTo$;%. System $S$ is said to be \textit{blocking} if it is not non--blocking;

\item \textit{accessible}, if for any $x\in X$ there exists a finite number of transitions
\[
x_{0} \rTo^{u_{1}} x_{1} \rTo^{u_{2}} \ldots \rTo^{u_{N}} x
\]
starting from an initial state $x_{0}$ in $X_{0}$ and ending up in $x$.

\end{itemize}

%For further use we denote by $\varnothing_{s}$ the empty system, so that the sets of states, of inputs and outputs are empty sets. 
We now introduce some notions which will be employed in the further developments. We start by introducing the notion of sub--system which formalizes the idea of extracting from the original system a subset of states, inputs and transitions.
%For such a transition $x\rTo^{u}x^{\prime }$, state $x^{\prime}$ is called a $u$-successor. Since \mbox{$\longrightarrow\subseteq X\times U\times X$} is a relation, for any state $x\in{X}$ and any input $u\in{U}$ there may be no $u$-successors, one, or many $u$-successors.

\begin{definition}
Given two systems $S_{1}=(X_{1},X_{0,1},U_{1},$ $\rTo_{1},Y_{1},H_{1})$ and $S_{2}=(X_{2},X_{0,2},U_{2},\rTo_{2},Y_{2},H_{2})$, system $S_{1}$ is a sub--system of $S_{2}$, denoted $S_{1} \sqsubseteq S_{2}$, if $X_{1}\subseteq X_{2}$, $X_{0,1}\subseteq X_{0,2}$, $U_{1}\subseteq U_{2}$, $\rTo_{1}\subseteq \rTo_{2}$, $Y_{1}\subseteq Y_{2}$ and $H_{1}(x)=H_{2}(x)$ for any $x\in X_{1}$.
\end{definition}

%The notion of sub--system naturally induces a preorder $\sqsubset$ on the class of systems so that $S_{1} \sqsubset S_{2}$ if $S_{1}$ is a sub--system of $S_{2}$. 
The following notion formalizes the idea of extracting the maximal non--blocking sub--system from a system, where maximality is given with respect to the notion of sub--system, which naturally induces a preorder on the class of systems.

\begin{definition}
Given a system $S=(X,X_{0},U,\rTo,Y,H)$ the non--blocking part of $S$ is a system $Nb(S)$ so that:
\begin{itemize}
\item[(i)] $Nb(S)$ is a non--blocking system;
\item[(ii)] $Nb(S)$ is a sub--system of $S$;
\item[(iii)] $S' \sqsubseteq Nb(S)$, for any non--blocking $S' \sqsubseteq S$.
\end{itemize}
\end{definition}

We finally introduce the notion of accessible part \cite{CassandrasBook} which formalizes the idea of extracting the maximal accessible sub--system from a system.

\begin{definition}
Given a system $S=(X,X_{0},U,\rTo,Y,H)$ the accessible part of $S$ is a system $Ac(S)$ so that:
\begin{itemize}
\item[(i)] $Ac(S)$ is an accessible system;
\item[(ii)] $Ac(S)$ is a sub--system of $S$;
\item[(iii)] $S' \sqsubseteq Ac(S)$, for any accessible $S' \sqsubseteq S$.
\end{itemize}

\end{definition}

%The maximal non--blocking sub--system of a system exists and it is unique.
In this paper we consider simulation and bisimulation relations \cite{Milner,Park} that are useful when analyzing or designing controllers for deterministic systems \cite{paulo}. Bisimulation relations are standard mechanisms to relate the properties of systems. Intuitively, a bisimulation relation
between a pair of systems $S_{1}$ and $S_{2}$ is a relation between the corresponding state sets explaining how a state trajectory $s_{1}$ of $S_{1}$ can be transformed into a state trajectory $s_{2}$ of $S_{2}$ and vice versa. While typical bisimulation relations require that $s_{1}$ and $s_{2}$ are
observationally indistinguishable, that is $H_{1}(s_{1}) = H_{2}(s_{2})$, we shall relax this by requiring $H_{1}(s_{1})$ to simply be close to $H_{2}(s_{2})$ where closeness is measured with respect to the metric on the output set. A simulation relation is a one-sided version of a bisimulation relation. The following notions have been introduced in \cite{AB-TAC07} and in a slightly different formulation in \cite{TabuadaTAC08}.

\begin{definition}
\label{ASR} Let \mbox{$S_{1}=(X_{1},X_{0,1},U_{1},\rTo_{1},Y_{1},H_{1})$} and \mbox{$S_{2}=(X_{2},X_{0,2},U_{2},\rTo_{2},Y_{2},H_{2})$} be metric systems
with the same output sets $Y_{1}=Y_{2}$ and metric $d$, and consider a precision $\varepsilon\in\mathbb{R}^{+}_{0}$. A relation
\[
\mathcal{R}\subseteq X_{1}\times X_{2},
\]
is said to be an $\varepsilon $--approximate simulation relation from $S_{1}$ to $S_{2}$, if the following conditions are satisfied:

\begin{itemize}
\item[(i)] for every $x_{1}\in X_{0,1}$, there exists $x_{2}\in X_{0,2}$ with
$(x_{1},x_{2})\in \mathcal{R}$;

\item[(ii)] for every $(x_{1},x_{2})\in \mathcal{R}$ we have \mbox{$d(H_{1}(x_{1}),H_{2}(x_{2}))\leq\varepsilon$};

\item[(iii)] for every $(x_{1},x_{2})\in \mathcal{R}$ we have that:

\mbox{$x_{1}\rTo_{1}^{u_{1}}x'_{1}$ in $S_{1}$} implies the existence of \mbox{$x_{2}\rTo_{2}^{u_{2}}x'_{2}$} in $S_{2}$ satisfying $(x^{\prime}%
_{1},x^{\prime}_{2})\in \mathcal{R}$.
\end{itemize}

System $S_{1}$ is \mbox{$\varepsilon$--approximately} simulated by $S_{2}$ or $S_{2}$ \mbox{$\varepsilon$--approximately} simulates $S_{1}$, denoted by
\mbox{$S_{1}\preceq_{\varepsilon}S_{2}$}, if there exists an \mbox{$\varepsilon$--approximate} simulation relation from $S_{1}$ to $S_{2}$.
When $\varepsilon=0$, system $S_{1}$ is said to be $0$--simulated by $S_{2}$ or
$S_{2}$ is said to $0$--simulate $S_{1}$.
\end{definition}

%The notion of $0$--simulation relations naturally induces a preorder $\sqsubset$ on the class of systems so that $S_{1} \sqsubseteq S_{2}$ if $S_{1} \preceq_{0} S_{2}$. 
By symmetrizing the notion of approximate simulation we obtain the notion of approximate bisimulation,
which is reported hereafter.

\begin{definition}
\label{BSR} Let \mbox{$S_{1}=(X_{1},X_{0,1},U_{1},\rTo_{1},Y_{1},H_{1})$} and \mbox{$S_{2}=(X_{2},X_{0,2},U_{2},\rTo_{2},Y_{2},H_{2})$} be metric systems with the same output sets $Y_{1}=Y_{2}$ and metric $d$, and consider a precision
$\varepsilon\in\mathbb{R}^{+}_{0}$. A relation
\[
\mathcal{R}\subseteq X_{1}\times X_{2},
\]
is said to be an \mbox{$\varepsilon$--approximate} bisimulation relation
between $S_{1}$ and $S_{2}$, if the following conditions are satisfied:

\begin{itemize}
\item[(i)] $\mathcal{R}$ is an \mbox{$\varepsilon$--approximate} simulation relation from
$S_{1}$ to $S_{2}$;

\item[(ii)] $\mathcal{R}^{-1}$ is an \mbox{$\varepsilon$--approximate} simulation
relation from $S_{2}$ to $S_{1}$.
\end{itemize}

System $S_{1}$ is \mbox{$\varepsilon$--approximately} bisimilar to $S_{2}$,
denoted by \mbox{$S_{1}\cong_{\varepsilon}S_{2}$}, if there exists an
$\varepsilon$--approximate bisimulation relation $\mathcal{R}$ between $S_{1}$ and $S_{2}$. When $\varepsilon=0$, system $S_{1}$ is said to be $0$--bisimilar or exactly bisimilar to $S_{2}$.
\end{definition}

We now introduce the notion of approximate composition of systems which is employed in the further developments to formalize the interconnection between a nonlinear control system representing the plant, and a symbolic system representing the symbolic controller.

\begin{definition}
\cite{TabuadaTAC08}\label{composition} Given two metric systems $S_{1}=(X_{1},X_{0,1},U_{1},\rTo_{1},Y_{1},H_{1})$ and $S_{2}=(X_{2},X_{0,2},U_{2},$ $\rTo_{2},Y_{2},H_{2})$, with the same output sets $Y_{1}=Y_{2}$ and metric $d$ and a precision $\varepsilon\in
\mathbb{R}_{0}^{+}$, the $\varepsilon$--approximate composition of $S_{1}$ and $S_{2}$ is the system:
\[
S_{1}\parallel_{\varepsilon}S_{2}:=(X,X_{0},U,\rTo,Y,H),
\]
where:

\begin{itemize}
\item $X=\{(x_{1},x_{2})\in X_{1}\times X_{2}:d(H_{1}(x_{1}),H_{2}(x_{2}))\leq \varepsilon\}$;

\item $X_{0}=X\cap(X_{0,1}\times X_{0,2})$;

\item $U=U_{1}\times U_{2}$;

\item $(x_{1},x_{2})\rTo^{(u_{1},u_{2})}(x_{1}^{\prime},x_{2}^{\prime})$ if $x_{1}\rTo_{1}^{u_{1}}x_{1}^{\prime}$ and $x_{2}\rTo_{2}^{u_{2}}x_{2}^{\prime}$;

\item $Y=Y_{1}$;

\item $H:X_{1}\times X_{2}\rightarrow Y$ is given by $H(x_{1},x_{2}):=H_{1}(x_{1})$, for any $(x_{1},x_{2})\in X$.
\end{itemize}
\end{definition}

The above notion of composition is asymmetric. This is because it models the interaction of systems $S_{1}$ and $S_{2}$ which play different roles in the composition. As it will be clarified in the next section, we interpret system $S_{1}$ as the plant system, i.e. the to--be--controlled process, and system $S_{2}$ as the controller.

\section{Problem Statement}\label{sec3}

In this paper we address the problem of symbolic control design for nonlinear systems with infinite states specifications modelled by differential equations. In order to formally define the control design problem under consideration, we first need to provide a formal notion of symbolic controllers. 
Given a control system \mbox{$\Sigma=(X,X_0,U,\mathcal{U},f)$} and a sampling time parameter $\tau\in{\mathbb{R}}^{+}$, we associate the following
system to $\Sigma$:
\begin{equation}
S_{\tau}(\Sigma):=(X,X_{0},\mathcal{U_{\tau}},\rTo_{\tau},Y,H),
\label{systemTD}
\end{equation}
where:

\begin{itemize}
\item $\mathcal{U_{\tau}}=\{u\in\mathcal{U}|$ the domain of
$u$ is $[0,\tau]$\};

\item $x\rTo_{\tau}^{u} x'$ if there exists
a trajectory $\xi:[0,\tau]\rightarrow X$ of $\Sigma$ satisfying \mbox{$\xi_{xu}(\tau)=x'$};

\item $Y=X$;

\item $H=1_{X}$.
\end{itemize}

System $S_{\tau}(\Sigma)$ is metric when we regard $Y=X$ as being equipped with the metric \mbox{$d(p,q)=\Vert p-q\Vert$}. The above system can be thought of as the time discretization of the control system $\Sigma$. 

\begin{definition}
Given the control system $\Sigma$, a sampling time $\tau\in\mathbb{R}^{+}$, a state quantization $\theta\in\mathbb{R}^{+}$ and an input quantization $\mu\in\mathbb{R}^{+}$, a \textit{symbolic controller} for $\Sigma$ is formalized by means of the system:
\[
C:=(X_{c},X_{c,0},U_{c},\rTo_{c},Y_{c},H_{c}),
\]
where\footnote{
The sets $[X]_{2\theta}$ and $[U]_{2\mu}$, are lattices embedded in the sets $\mathbb{R}^{n}$ and $U$, with precisions $\theta$ and $\mu$ respectively, as formally defined in the Appendix.}:

\begin{itemize}
\item $X_{c}=[X]_{2\theta}$;

\item $X_{c,0}\subseteq X_{c}$;

\item $U_{c}=\{u\in\mathcal{U}_{\tau}|$ the co--domain of $u$
is $[U]_{2\mu}$\};

\item $\rTo_{c} \subseteq X_{c} \times U_{c} \times X_{c}$;

\item $Y_{c}=X_{c}$;

\item $H_{c}=1_{X_{c}}$.
\end{itemize}

\end{definition}

\begin{figure}
\begin{center}
\includegraphics[scale=0.6]{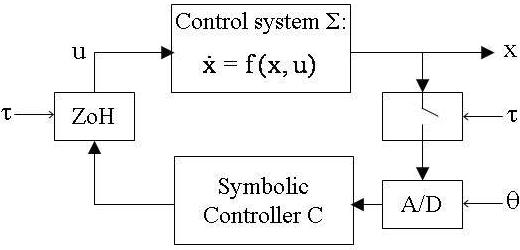}
\caption{Approximate composition of the plant and specification systems.} 
\label{fig0}
\end{center}
\end{figure}

We denote by $\mathcal{C}^{\tau,\theta,\mu}(\Sigma)$ the class of symbolic controllers with sampling time $\tau$, state quantization $\theta$ and input  quantization $\mu$, associated with $\Sigma$. 
The $\theta$--approximate composition between the time discretization $S_{\tau}(\Sigma)$ of a control system $\Sigma$ and a symbolic controller $C\in\mathcal{C}^{\tau,\theta,\mu}(\Sigma)$ formalizes classical static state feedback control schemes with digital controllers, studied in the literature, see e.g. \cite{DigitalBook}, as illustrated in Figure \ref{fig0}: The state signal $\xi_{x_{0}u}(t)$ at time $t\in\mathbb{R}^{+}$ is firstly sampled with sampling time $\tau\in\mathbb{R}^{+}$, then quantized through an Analog--to--Digital (A/D) converter with precision $\theta\in\mathbb{R}^{+}$ which associates to a state $\xi_{x_{0}u}(\tau)$, the unique state $x\in X_{c}$ for which\footnote{The set $\mathcal{B}_{[\theta[}(x)$ is defined in the Appendix.} $\xi_{x_{0}u}(\tau)\in\mathcal{B}_{[\theta[}(x)$; the obtained digital/symbolic signal is then plugged as input to the digital/symbolic controller $C$ which outputs a symbolic signal taking values in $[U]_{2\mu}$. Such symbolic signal is then plugged into a Zero order Holder (ZoH) with sampling time parameter $\tau$ which outputs in turn, a piecewise--constant signal $u$ that is finally plugged as digital/symbolic control input to the control system $\Sigma$.
%When the state $x$ of the control system $\Sigma$ is in\footnote{The set $\mathcal{B}_{[\theta)}(y)$ is defined in the Appendix.} $\mathcal{B}_{[\theta)}(y)$ for some $y\in X_{c}$ at time $t=k\tau$ for $k\in\mathbb{N}$, the symbolic controller $C$ provides a control input $u\in U_{c}$ so that $y \rTo_{c}^{u} z$ for some $z \in X_{c}$. The interaction between the control system $\Sigma$ and the symbolic controller $C$ is formally captured by the approximate composition $S_{\tau}(\Sigma)\parallel_{\theta} C$.
%%More precisely, the feedback composition of $\Sigma$ and $S_{c}$ is given by $S_{\tau}(\Sigma)\parallel_{\theta} S_{c}$.

We are now ready to formally state the symbolic control design problem that we consider in this paper. Consider a plant nonlinear control system:
\begin{equation}
P=(X_{p},X_{p,0},U_{p},\mathcal{U}_{p},f_{p}),
\label{plant}
\end{equation}
and a specification nonlinear autonomous system:
\[
Q=(X_{q},X_{q,0},g_{q}).
\]
For the sake of homogeneity in the notation of the plant $P$ and the specification $Q$ we rephrase the above tuple by means of:
\begin{equation}
Q=(X_{q},X_{q,0},U_{q},\mathcal{U}_{q},f_{q}),
\label{specif}
\end{equation}
where $U_{q}=\{u_{q}\}$ with $u_{q}=0$, $\mathcal{U}_{q}=\{\mathbf{u}_{q}\}$ with $\mathbf{u}_{q}=\mathbf{0}$, the signal $\mathbf{0}$ being the identically null function, and $f_{q}(x,u)=g_{q}(x)+u$ for any $(x,u)\in X_{q}\times U_{q}$.

\begin{problem}
\label{problem}
Given a plant nonlinear control system $P$ as in (\ref{plant}), a specification nonlinear autonomous system $Q$ as in (\ref{specif}) and a desired precision $\varepsilon\in\mathbb{R}^{+}$, find quantization parameters $\tau,\theta,\mu\in\mathbb{R}^{+}$ and a symbolic controller $C\in\mathcal{C}^{\tau,\theta,\mu}(P)$ such that:

\begin{itemize}
\item[(i)] $(S_{\tau}(P)\parallel_{\theta}C) \preceq_{\varepsilon} S_{\tau}(Q)$;

\item[(ii)] $S_{\tau}(P)\parallel_{\theta}C$ is non--blocking.

\end{itemize}
\end{problem}

The above control design problem asks for a symbolic controller $C$ that implements the behaviour of the specification $Q$, up to a precision $\varepsilon$ that can be chosen as small as desired. In other words, in Problem \ref{problem} we look for a symbolic controller $C$ so that the approximate composition between the plant $P$ and the controller $C$ satisfies or conforms \cite{ModelChecking} the specification $Q$ with an arbitrarily small precision. The symbolic controller is furthermore requested to be non--blocking in order to prevent occurrence of deadlocks in the interaction between the plant and the symbolic controller. This control design problem can be seen as an approximated version of similarity games, as discussed in \cite{paulo}. Similar problems have been studied in the literature (in a non--approximating settings) in the context of supervisory control \cite{CassandrasBook}, symbolic control design for piecewise--affine systems enforcing temporal logic specifications \cite{BeltaCDC09}, among many other work.

\section{Symbolic Control Design with Infinite States Specifications}\label{sec4}
In this section we provide the solution to Problem \ref{problem}. Inspired by the so--called correct--by design approach, see e.g. \cite{LTLControl,TabuadaTAC08,BeltaCDC09}, we first construct the symbolic systems associated with the plant $P$ and the specification $Q$ in Section \ref{subsection1}, we then solve the control design problem at the symbolic layer in Section \ref{subsection2} to finally come back at the continuous layer in Section \ref{subsection3} by providing the bounds in the approximation scheme that we propose, which guarantee the solution to Problem \ref{problem}.

\subsection{From the Continuous Layer to the Symbolic Layer}\label{subsection1}
In this section we present some results based on the work of \cite{PolaAutom2008} for constructing symbolic systems associated with the plant $P$ and the specification $Q$. We start by recalling from \cite{IncrementalS}, the notion of incremental input--to--state stability for nonlinear control systems.
%
%\begin{definition}
%A control system $\Sigma$ is incrementally globally asymptotically stable ($\delta$-GAS) if it is forward complete and there exists a $\mathcal{KL}$ function $\beta$ such that for any $t\in{\mathbb{R}_0^+}$, any $x,x'\in{\mathbb{R}^n}$ and any $\upsilon\in\mathcal{U}$ the following condition is satisfied:
%\begin{equation}
%\left\Vert\xi_{x\upsilon}(t)-\xi_{x'\upsilon}(t)\right\Vert \leq\beta\left(\left\Vert x-x'\right\Vert,t\right). \label{delta_GAS}%
%\end{equation}
%\end{definition}

\begin{definition}
\label{dISS}
A control system $\Sigma$ is incrementally input--to--state stable ($\delta$--ISS) if it is forward complete and there exist a $\mathcal{KL}$ function $\beta$ and a $\mathcal{K}_{\infty}$ function $\gamma$ such that for any $t\in{\mathbb{R}_0^+}$, any $x,x'\in{\mathbb{R}^n}$, and any $u$, $u'\in\mathcal{U}$ the following condition is satisfied:
\begin{equation}
\left\Vert \xi_{xu}(t)-\xi_{x'u'}(t)\right\Vert \leq\beta\left(\left\Vert x-x'\right\Vert,t\right)+\gamma\left(\left\Vert u-u'\right\Vert_{\infty}\right). \label{delta_ISS}%
\end{equation}
\end{definition}

A characterization of the above incremental stability notion in terms of dissipation inequalities can be found in \cite{IncrementalS}.
%It is readily seen from (\ref{delta_GAS}) and (\ref{delta_ISS}) that $\delta$-ISS implies $\delta$-GAS while the converse is not true in general. Moreover for autonomous nonlinear systems $\delta$-ISS and $\delta$-GAS are equivalent.\\
Given a $\delta$--ISS nonlinear control system $\Sigma$ of the form (\ref{NLCsystem}), a sampling time $\tau\in\mathbb{R}^{+}$, a state quantization $\eta\in\mathbb{R}^{+}$ and an input quantization $\mu\in\mathbb{R}^{+}$ consider the following system:
\begin{equation}
S_{\tau,\eta,\mu}(\Sigma):=(X_{\tau,\eta,\mu},X_{0,\tau,\eta,\mu},U_{\tau,\eta,\mu},\rTo_{\tau,\eta,\mu},Y_{\tau,\eta,\mu},H_{\tau,\eta,\mu}),
\label{symbmodel}
\end{equation}
where:

\begin{itemize}
\item $X_{\tau,\eta,\mu}=[X]_{2\eta}$;
\item $X_{0,\tau,\eta,\mu}=X_{\tau,\eta,\mu}\cap X_{0}$;
\item $U_{\tau,\eta,\mu}=[U]_{2\mu}$;
\item $x\rTo_{\tau,\eta,\mu}^{u} y$ if $\xi_{xu}(\tau)\in\mathcal{B}_{[\eta[}(y)\cap X$;
\item $Y_{\tau,\eta,\mu}=X$;
\item $H_{\tau,\eta,\mu}=\imath : X_{\tau,\eta,\mu}\hookrightarrow Y_{\tau,\eta,\mu}$.
\end{itemize}

It is readily seen from the definition of $X_{\tau,\eta,\mu}$ and $U_{\tau,\eta,\mu}$ that system $S_{\tau,\eta,\mu}(\Sigma)$ is countable and becomes symbolic when the state space $X$ and the input space $U$ are bounded sets. System $S_{\tau,\eta,\mu}(\Sigma)$ is basically equivalent to the symbolic model proposed in \cite{PolaAutom2008}. The main difference is that, while the symbolic model in \cite{PolaAutom2008} is not guaranteed to be deterministic, system $S_{\tau,\eta,\mu}(\Sigma)$ is so, as formally stated in the following result:

\begin{proposition}
\label{prop:determinism}
System $S_{\tau,\eta,\mu}(\Sigma)$ is deterministic.
\end{proposition}

\begin{proof}
The existence and uniqueness of a trajectory from an initial condition $x\in X_{\tau,\eta,\mu}$ with input $u \in U_{\tau,\eta,\mu}$ guarantees that
$\xi_{xu}(\tau)$ is uniquely determined. Since the collection of sets $\{\mathcal{B}_{[\eta[}(y)\cap X\}_{y\in X_{\tau,\eta,\mu}}$ is a partition of  $X$, there exists at most one state $y\in X_{\tau,\eta,\mu}$\ such that $\xi_{xu}(\tau)\in\mathcal{B}_{[\eta[}(y)\cap X$.
\end{proof}

We stress that determinism in the symbolic system $S_{\tau,\eta,\mu}(\Sigma)$ is an important property because algorithmic synthesis of symbolic systems simplifies when systems are deterministic \cite{paulo}. %Moreover, determinism in $S_{\tau,\eta,\mu}(\Sigma)$ reduces a number of spurious transitions which are instead included in the symbolic systems in \cite{PolaAutom2008} and which are not essential in control design.???????\\
We can now give the following result that establishes sufficient conditions for the existence and construction of symbolic systems for nonlinear control systems.

\begin{theorem}
Consider a $\delta$--ISS nonlinear control system $\Sigma=(X,X_{0},U,\mathcal{U},f)$ and a desired precision \mbox{$\theta\in\mathbb{R}^{+}$}. For any sampling time $\tau\in\mathbb{R}^{+}$, state quantization $\eta\in\mathbb{R}^{+}$ and input quantization $\mu\in\mathbb{R}^{+}$ satisfying the following inequality:
\begin{equation}
\beta(\theta,\tau)+\gamma(\mu)+\eta \leq \theta ,
\label{condAUT}
\end{equation}
systems $S_{\tau,\eta,\mu}(\Sigma)$ and $S_{\tau}(\Sigma)$ are $\theta$--approximately bisimilar.
\label{polaut}
\end{theorem}

\begin{proof}
The proof of the above result can be given along the lines of Theorem 5.1 in \cite{PolaAutom2008}. We include it here for the sake of completeness.
Consider the relation $\mathcal{R}\subseteq X\times X_{\tau,\eta,\mu}$ defined by $(x,y)\in \mathcal{R}$ if and only if $x\in\mathcal{B}_{[\theta[}(y)\cap X$. We start by showing that condition (i) of Definition \ref{ASR} holds. Consider an initial condition $x_{0}\in X_{0}$. By definition of the set
$X_{0,\tau,\eta,\mu}$ there exists $y_{0}\in X_{0,\tau,\eta,\mu}$ so that $(x_{0},y_{0}) \in \mathcal{R}$.
Condition (ii) in Definition \ref{ASR} is satisfied by the definition of $\mathcal{R}$. Let us now show that condition (iii) in Definition \ref{ASR} holds. Consider any $(x,y)\in \mathcal{R}$. Consider any $u_{1}\in \mathcal{U}_{\tau}$ and the transition $x\rTo_{\tau}^{u_{1}} w$ in $S_{\tau}(\Sigma)$. There exists $u_{2}\in U_{\tau,\eta,\mu}$ such that:
\begin{equation}
\Vert u_{2} - u_{1} \Vert_{\infty} \leq \mu. \label{a5}%
\end{equation}
Set $z=\xi_{yu_{2}}(\tau)$. Since $X={\textstyle\bigcup \nolimits_{v\in X_{\tau,\eta,\mu}}}\mathcal{B}_{[\eta[}(v) \cap X$, there exists $v\in X_{\tau,\eta,\mu}$ such that:
\begin{equation}
z\in\mathcal{B}_{[\eta[}(v), \label{a2}%
\end{equation}
and therefore $y \rTo_{\tau,\eta,\mu}^{u_{2}} v$ in $S_{\tau,\eta,\mu}(\Sigma)$. Since
$\Sigma$ is $\delta$--ISS, by the definition of $\mathcal{R}$ and by condition (\ref{a5}), the following chain of inequalities holds:
\[
\Vert w-z \Vert\leq\beta(\Vert x-y\Vert,\tau)+\gamma(\Vert u_{1}-u_{2}\Vert_{\infty})\leq\beta(\theta,\tau)+\gamma(\mu),
\]
which implies:
\begin{equation}
w\in
%\operatorname*{cl}\left(  \mathcal{B}_{\beta(\Vert x-y\Vert,\tau )+\gamma(\Vert u_{1}-u_{2}\Vert_{\infty})}(z)\right)  \subseteq
\mathcal{B}_{\beta(\theta,\tau)+\gamma(\mu)}(z).\label{a6}%
\end{equation}
By combining the inclusions in (\ref{a2}) and (\ref{a6}), it is readily seen that $w\in\mathcal{B}_{[\beta(\theta,\tau)+\gamma(\mu)+\eta[}(v)$.
By the inequality in (\ref{polaut}), $\mathcal{B}_{[\beta(\theta,\tau)+\gamma(\mu)+\eta[}(v)\subseteq\mathcal{B}_{[\theta[}(v)$, which implies $(w,v)\in \mathcal{R}$ and hence, condition (iii) in Definition \ref{ASR} holds. Thus, condition (i) in Definition \ref{BSR} is satisfied. By using similar arguments it is possible to show condition (ii) of Definition \ref{BSR}.
\end{proof}

The above result is conceptually equivalent to Theorem 5.1 in \cite{PolaAutom2008}. The main difference is that while Theorem \ref{polaut} relates nonlinear systems to \textit{deterministic} symbolic systems, Theorem 5.1 in \cite{PolaAutom2008} relates nonlinear systems to symbolic models which are in general nondeterministic.\\
The above result is now employed to define symbolic systems for the plant and the specification. Consider a plant system $P$ as defined in (\ref{plant}) and a specification system $Q$ as defined in (\ref{specif}). Suppose that $P$ and $Q$ are $\delta$--ISS and choose a precision $\theta_{p}\in\mathbb{R}^{+}$ and a precision $\theta_{q}\in\mathbb{R}^{+}$, required in the construction of the symbolic systems for $P$ and $Q$, respectively. Let $\beta_{p}$ and $\gamma_{p}$ be a $\mathcal{KL}$ function and a $\mathcal{K}_{\infty}$ function guaranteeing the $\delta$--ISS stability property for $P$ and $\beta_{q}$ be a $\mathcal{KL}$ function guaranteeing the $\delta$--ISS stability property for $Q$. Find quantization parameters $\tau,\eta,\mu\in\mathbb{R}^{+}$ such that:
\begin{eqnarray}
& & \beta_{p}(\theta_{p},\tau)+\gamma_{p}(\mu)+\eta\leq \theta_{p},\nonumber\\
& & \beta_{q}(\theta_{q},\tau)+\eta\leq \theta_{q},
\label{condAUT2}
\end{eqnarray}

It is readily seen that parameters $\tau,\eta,\mu\in\mathbb{R}^{+}$ satisfying the above inequalities always exist. By Theorem \ref{polaut}, $S_{\tau,\eta,\mu}(P)$ is $\theta_{p}$--approximately bisimilar to $S_{\tau}(P)$ and $S_{\tau,\eta,0}(Q)$ is $\theta_{q}$--approximately bisimilar to $S_{\tau}(Q)$. For the sake of notational simplicity in the further developments we refer to the systems $S_{\tau,\eta,\mu}(P)$ and $S_{\tau,\eta,0}(Q)$, by means of $S_{p}$ and $S_{q}$, respectively.

\subsection{Control Design at the Symbolic Layer}\label{subsection2}

Problem \ref{problem} translates to the following problem at the symbolic layer:

\begin{problem}
\label{problemD} Given system $S_{p}$ and system $S_{q}$, find a symbolic controller $C\in \mathcal{C}^{\tau,\theta,\mu}(P)$ such that:

\begin{itemize}
\item[(i)] $(S_{p}\parallel_{0}C) \preceq_{0} S_{q}$;

\item[(ii)] $S_{p}\parallel_{0}C$ is non--blocking.

\end{itemize}
\end{problem}

We start by introducing a technical lemma that will be used in the sequel.

\begin{lemma}
\label{lemma}
Consider three metric systems $S_{i}=(X_{i},X_{0,i},U_{i},\rTo_{i},Y,H_{i})$, $i=1,2,3$. The following properties hold:

\begin{itemize}
\item[(i)] \cite{AB-TAC07} For all $\varepsilon_{1}\in\mathbb{R}^{+}_{0}$, if
$S_{1}\preceq_{\varepsilon_{1}}S_{2}$ then $S_{1}\preceq_{\varepsilon_{2}}S_{2}$, for all $\varepsilon_{2}\geq\varepsilon_{1}$;

\item[(ii)] \cite{AB-TAC07} For all $\varepsilon_{1},\varepsilon
_{2}\in\mathbb{R}_{0}^{+}$, if $S_{1}\preceq_{\varepsilon_{1}}S_{2}$ and
$S_{2}\preceq_{\varepsilon_{2}}S_{3}$, then $S_{1}\preceq_{\varepsilon
_{1}+\varepsilon_{2}}S_{3}$;

\item[(iii)] For all $\varepsilon\in\mathbb{R}^{+}_{0}$, $S_{1}\parallel_{\varepsilon}S_{2}%
\preceq_{\varepsilon} S_{2}$.

\end{itemize}
\end{lemma}

\begin{proof}
[Proof of (iii)] Denote $S_{1}\parallel_{\varepsilon}S_{2}$ by the tuple $(X,X_{0},U,\rTo,Y,H)$ and define: 
\[
\mathcal{R=}\{((x_{1},x_{2}),x)\in X\times X_{2}:x_{2}=x\}.
\]
We start by showing that condition (i) in Definition \ref{ASR} holds. Consider any initial condition $(x_{0,1},x_{0,2})\in X_{0}$. Since $x_{0,2}\in X_{2}$, by choosing $x_{0}=x_{0,2}$ we have that $((x_{0,1},x_{0,2}),x_{0})\in\mathcal{R}$.
We now show that also condition (ii) in Definition \ref{ASR} holds. Consider any $((x_{1},x_{2}),x)\in\mathcal{R}$. Since $x_{2}=x$, then $H_{2}(x_{2})=H_{2}(x)$, hence by Definition \ref{composition} of approximate composition $d(H(x_{1},x_{2}),H_{2}(x))=d(H_{1}(x_{1}),H_{2}(x_{2}))\leq\varepsilon$.
We conclude by showing that condition (iii) in Definition \ref{ASR} holds.
Consider any $((x_{1},x_{2}),x)\in\mathcal{R}$ and any transition $(x_{1},x_{2})\rTo^{(u_{1},u_{2})}(x_{1}^{\prime},x_{2}%
^{\prime})$ in $S_{1}\parallel_{\varepsilon} S_{2}$.
%By definition of composition, there exist a pair of transitions $x_{1}\rTo^{u_{1}}x_{1}^{\prime}$ in $S_{1}$ and $x_{2}\rTo^{u_{2}}x_{2}^{\prime}$ in $S_{2}$ so that $d(H_{1}(x_{1}^{\prime}),H_{2}(x_{2}^{\prime}))\leq\varepsilon$.
Choose the transition $x\rTo^{u_{2}}_{2}x^{\prime}$ in $S_{2}$ so that $x^{\prime}=x_{2}^{\prime}$. By definition of the systems involved such transition exists. This implies that $((x_{1}^{\prime},x_{2}^{\prime}),x^{\prime})\in\mathcal{R}$, which concludes the proof.
\end{proof}

We are now ready to provide the solution to Problem \ref{problemD}. Define:
\begin{equation}
C^{\ast}=S_{p}\parallel_{0}S_{q}.
\label{controller}
\end{equation}

\begin{theorem}
$Nb(C^{\ast})$ solves Problem \ref{problemD}.
\end{theorem}

\begin{proof}
We start by proving condition (i) of Problem \ref{problemD}. By Lemma \ref{lemma} (iii), we obtain:
\begin{equation}
S_{p}\parallel_{0}Nb(C^{\ast})\preceq_{0} Nb(C^{\ast}).
\label{cond111}
\end{equation}
By the definition of $Nb(C^{\ast})$ it is readily seen that:
\begin{equation}
Nb(C^{\ast})\preceq_{0} C^{\ast}.
\label{cond222}
\end{equation}
By the definition of $C^{\ast}=S_{p}\parallel_{0}S_{q}$ and Lemma \ref{lemma} (iii), one gets:
\begin{equation}
C^{\ast}\preceq_{0}S_{q}.
\label{cond333}
\end{equation}
By combining conditions in (\ref{cond111}), (\ref{cond222}), (\ref{cond333}) and by Lemma \ref{lemma} (ii) we obtain:
\[
S_{p}\parallel_{0}Nb(C^{\ast})\preceq_{0} S_{q}.
\]
Hence, condition (i) of Problem \ref{problemD} is proved. We now prove condition (ii) of Problem \ref{problemD}. Consider any state $(p_{1},p_{2},q)$ of $S_{p}\parallel_{0}Nb(C^{\ast})$. Since $Nb(C^{\ast})$ is non--blocking there exists a state $(p_{2}^{+},q^{+})$ of $Nb(C^{\ast})$ so that $(p_{2},q)\rTo^{u}(p_{2}^{+},q^{+})$ is a transition of $Nb(C^{\ast})$ for some input $u=(u_{2},u_{3})$. Since $Nb(C^{\ast})$ is a sub--system of $C^{\ast}=S_{p}\parallel_{0}S_{q}$, then by choosing $p_{1}^{+}=p_{2}^{+}$ and $u_{1}=u_{2}$, the transition $p_{1} \rTo^{u_{1}} p_{1}^{+}$ is a transition of $S_{p}$. Since by construction $p_{1}^{+}=p_{2}^{+}$ then
$(p_{1}^{+},p_{2}^{+},q^{+})$ is a state of $S_{p}\parallel_{0}Nb(C^{\ast})$ and therefore $(p_{1},p_{2},q)\rTo^{(u_{1},u)} (p_{1}^{+},p_{2}^{+},q^{+})$ is a transition of $S_{p}\parallel_{0}Nb(C^{\ast})$, which concludes the proof.
\end{proof}

We conclude this section by showing that the controller $Nb(C^{\ast})$ is the maximal system solving Problem \ref{problemD} in the sense of the preorder  naturally induced by the notion of $0$--simulation relations.

\begin{theorem}
For any system $C$ solving Problem \ref{problemD}
\[
(S_{p}\parallel_{0}C)\preceq_{0}(S_{p}\parallel_{0}Nb(C^{\ast})).
\]
\end{theorem}

\begin{proof}
Denote by $S_{p}^{1}$ and $S_{p}^{2}$ copies of $S_{p}$ that are connected to $C$ and $Nb(C^{\ast})$, respectively; denote by $X_{pc}$ and $X_{pc^{\ast}}$ the state spaces of $S_{p}^{1}\parallel_{0}C$ and $S_{p}^{2}\parallel_{0}Nb(C^{\ast})$ and by $X_{pc}^{0}$ and $X_{pc^{\ast}}^{0}$ the corresponding sets of initial states. Moreover let $C^{\ast}=S_{p}^{c}\parallel_{0}S_{q}^{c}$, where
$S_{p}^{c}$ and $S_{q}^{c}$ are the copies of $S_{p}$ and $S_{q}$ in the
controller and define:
\[
\mathcal{R}=
\{((p_{1},c),(p_{2},p_{3},q))  \in X_{pc}\times X_{pc^{\ast}}:
((p_{1},c),q) \in \mathcal{R}_{1} \wedge p_{1}=p_{2}\},
\]
where $\mathcal{R}_{1}$ is a $0$--simulation relation from $S_{p}^{1} \parallel_{0} C$ to $S_{q}$. We start by showing that condition (i) in Definition \ref{ASR} holds. Consider any initial condition $(p_{1}^{0},c^{0}) \in X_{pc}^{0}$. Since $(S_{p}\parallel_{0}C)\preceq_{0}S_{q}$ there exists $q^{0}\in S_{q}$ s.t. $((p_{1}^{0},c^{0}),q^{0})\in \mathcal{R}_{1}$. By choosing $p_{2}^{0}=p_{3}^{0}=p_{1}^{0}$, we have $(p_{2}^{0},p_{3}^{0},q^{0})\in X_{pc^{\ast}}^{0}$ and hence, $((p_{1}^{0},c^{0}),(p_{2}^{0},p_{3}^{0},q^{0}))\in\mathcal{R}$. We now show that also condition (ii) in Definition \ref{ASR} holds.
%Consider
%any $((p_{1},c),(p_{2},p_{3},q))\in\mathcal{R}$ and let $H_{pc}$ and $H_{ppq}$ be the output functions of
%$S_{p}^{1}\parallel_{0}C$ and $S_{p}^{2}\parallel_{0}Nb(C^{\ast})$, respectively.
%Since $((p_{1},c),q)  \in\mathcal{R}_{1}$, then
%$H_{pc}(p_{1},c)=H_{p}(p_{1})=H_{q}(q)$;
%moreover, by definition of exact composition, $H_{ppq}(p_{2},p_{3},q)=H_{p}(p_{2})=H_{p}(p_{3})
%=H_{q}(q)$.
Since $H_{p}(p_{1})=H_{p}(p_{2})$, we can conclude $d(H_{pc}(p_{1},c),H_{ppq}(p_{2},p_{3},q))=d
(H_{p}(p_{1}),H_{p}(p_{2}))=0$.
We conclude by showing that condition (iii) in Definition \ref{ASR} holds.
Consider any $((p_{1},c),(p_{2},p_{3},q))  \in\mathcal{R}$ and any transition
$(p_{1},c)\rTo^{(u_{p},u_{c})}(p_{1}^{+},c^{+})$ in $S_{p}^{1}\parallel_{0}C$.
%By definition of composition, there exists a pair of transitions $p_{1}\rTo^{u_{p}}p_{1}^{+}$ and $c\rTo^{u_{c}}c^{+}$ in $S_{p}^{1}$ and $C$, respectively.
Since $((p_{1},c),q)\in\mathcal{R}_{1}$,
there exists a transition $q\rTo^{v}q^{+}$ in $S_{q}$ so that
$((p_{1}^{+},c^{+}),q^{+})\in\mathcal{R}_{1}$.
Hence $H_{pc}(p_{1}^{+},c^{+})=H_{p}(p_{1}^{+})=H_{q}(q^{+})$ and $q^{+}=p_{1}^{+}$.
Now, since $p_{1}=p_{2}=p_{3}=q$, we consider the transitions
$p_{2}\rTo^{u_{p}}p_{2}^{+}$ in $S_{p}^{2}$, $p_{3}\rTo^{u_{p}}p_{3}^{+}$ in $S_{p}^{c}$ and $q\rTo^{v}q^{+}$ in $S_{q}^{c}$ with $p_{3}^{+}=p_{2}^{+}=p_{1}^{+}=q^{+}$. Notice that such
transitions exist. Hence $(p_{2}^{+},p_{c}^{+},q^{+})$ is a state of $S_{p}^{2}\parallel_{0}Nb(C^{\ast})$ and the transition
$(p_{2},p_{3},q) \rTo^{(u_{p},u_{p},v)}(p_{2}^{+},p_{3}^{+},q^{+})$ is in
$S_{p}^{2}\parallel_{0} Nb(C^{\ast})$, which implies $((p_{1}^{+},c^{+}),(p_{2}^{+},p_{3}^{+},q^{+}))\in\mathcal{R}$.
\end{proof}

This result is important because it shows that the controller $Nb(C^{\ast})$ implements the maximal non--blocking behaviour of the specification symbolic system $S_{q}$, which can be implemented by the plant symbolic system $S_{p}$.

\subsection{From the Symbolic Layer to the Continuous Layer}\label{subsection3}
We now have all the ingredients to present one of the main results of this paper which shows that there exists an appropriate choice of quantization parameters so that the symbolic controller $Nb(C^{\ast})$ with $C^{\ast}$ defined in (\ref{controller}) solves Problem \ref{problem}.

\begin{theorem}
Consider the plant system $P$ as in (\ref{plant}), the specification system $Q$ as in (\ref{specif}) and a precision $\varepsilon\in\mathbb{R}^{+}$.
Suppose that $P$ and $Q$ are $\delta$--ISS and choose parameters $\theta_{p},\theta_{q}\in\mathbb{R}^{+}$ so that:
\begin{equation}
\theta_{p}+\theta_{q}\leq\varepsilon.
\label{cond1111}
\end{equation}
Furthermore choose parameters $\tau,\eta,\mu\in\mathbb{R}^{+}$ satisfying the inequalities in (\ref{condAUT2}). Then the symbolic controller $Nb(C^{\ast})\in\mathcal{C}^{\tau,\theta,\mu}(P)$ with $\theta=\theta_{p}$ and $C^{\ast}$ defined in (\ref{controller}) with $S_{p}=S_{\tau,\eta,\mu}(P)$ and $S_{q}=S_{\tau,\eta,0}(Q)$, solves Problem \ref{problem}.
\label{ThMain1}
\end{theorem}

\begin{proof}
We start by proving condition (i) of Problem \ref{problem}. By Lemma \ref{lemma} (iii), we obtain:
\begin{equation}
S_{\tau}(P)\parallel_{\theta_{p}}Nb(C^{\ast})\preceq_{\theta_{p}} Nb(C^{\ast}).
\label{cond11}
\end{equation}
By the definition of $Nb(C^{\ast})$ it is readily seen that:
\begin{equation}
Nb(C^{\ast})\preceq_{0} C^{\ast}.
\label{cond22}
\end{equation}
By the definition of $C^{\ast}=S_{p}\parallel_{0}S_{q}$ and Lemma \ref{lemma} (iii), one gets:
\begin{equation}
C^{\ast}\preceq_{0}S_{q}.
\label{cond33}
\end{equation}
Since $S_{q}$ is $\theta_{q}$--approximately bisimilar to $S_{\tau}(Q)$ then:
\begin{equation}
S_{q} \preceq_{\theta_{q}} S_{\tau}(Q).
\label{cond44}
\end{equation}
By combining conditions in (\ref{cond11}), (\ref{cond22}), (\ref{cond33}), (\ref{cond44}) and by Lemma \ref{lemma} (ii) we obtain:
\[
S_{\tau}(P)\parallel_{\theta_{p}}Nb(C^{\ast})\preceq_{\theta_{p}+\theta_{q}} S_{\tau}(Q).
\]
Since by (\ref{cond1111}), $\theta_{p}+\theta_{q}\leq \varepsilon$, by Lemma \ref{lemma} (i), condition (i) of Problem \ref{problem} is proved.
We now prove condition (ii) of Problem \ref{problem}. Consider any state $(p_{1},p_{2},q)$ of $S_{\tau}(P)\parallel_{\theta_{p}}Nb(C^{\ast})$. Since $Nb(C^{\ast})$ is non--blocking there exists a state $(p_{2}^{+},q^{+})$ of $Nb(C^{\ast})$ so that $(p_{2},q)\rTo^{u}(p_{2}^{+},q^{+})$ is a transition of $Nb(C^{\ast})$ for some input $u=(u_{2},u_{3})$. Since $S_{\tau}(P)$ and $S_{p}$ are $\theta_{p}$--approximately bisimilar, for the transition $p_{2}\rTo^{u_{2}} p_{2}^{+}$ in $S_{p}$ there exists a transition $p_{1}\rTo^{u_{1}} p_{1}^{+}$ in $S_{\tau}(P)$ so that $d(H_{p}(p_{1}^{+}),H_{p}(p_{2}^{+}))\leq \theta_{p}$. This implies that $(p_{1}^{+},p_{2}^{+},q^{+})$ is a state of $S_{\tau}(P)\parallel_{\theta_{p}}Nb(C^{\ast})$ and therefore that $(p_{1},p_{2},q)\rTo^{(u_{1},u)}(p_{1}^{+},p_{2}^{+},q^{+})$ is a transition of $S_{\tau}(P)\parallel_{\theta_{p}}Nb(C^{\ast})$, which concludes the proof.
\end{proof}

\section{Integrated Symbolic Control Design}{\label{sec5}

The construction of the symbolic controller $Nb(S_{p}\parallel_{0}S_{q})$ solving Problem \ref{problem} relies upon the basic--steps procedure illustrated in Algorithm \ref{alg3}.

\incmargin{1em}
\restylealgo{boxed}\linesnumbered
\begin{algorithm}[H]
\label{alg3}
\SetLine
\caption{Construction of $Nb(S_{p}\parallel_{0}S_{q})$.}
Construct system $S_{p}$, $\theta_{p}$--approximately bisimilar to $S_{\tau}(P)$\;
Construct system $S_{q}$, $\theta_{q}$--approximately bisimilar to $S_{\tau}(Q)$\;
Construct the composition $S_{p}\parallel_{0}S_{q}$\;
Compute the non--blocking part $Nb(S_{p}\parallel_{0}S_{q})$ of $S_{p}\parallel_{0}S_{q}$.
\end{algorithm}
\decmargin{1em}

%The above procedure which shares similar ideas with other approaches currently available in the literature on symbolic control design of continuous and hybrid systems, see e.g. \cite{LTLControl,BeltaCDC09,TabuadaTAC08}, is not efficient from the computational point of view, as discussed hereafter. 

The procedure in Algorithm \ref{alg3} is common with other approaches currently available in the literature for symbolic control design of continuous and hybrid systems, see e.g. \cite{LTLControl,BeltaCDC09,TabuadaTAC08}. Software implementation of Algorithm \ref{alg3} requires that:

\begin{itemize}
\item State space $X_{p}$ and set of input values $U_{p}$ of $P$ are bounded;
\item State space $X_{q}$ of $Q$ is bounded.
\end{itemize}

The above assumptions, while being reasonable in many realistic engineering control problems, are also needed to store the transitions of systems $S_{p}$ and $S_{q}$ in a computer machine, whose memory resources are limited by their nature. In this section, we suppose that the plant $P$ and the specification $Q$ satisfy the above assumptions.
%The construction of $S_{p}$ is executed in $O(card([X_{p}]_{2\eta})\cdot card([U_{p}]_{2\mu}))$. The construction of $S_{q}$ is executed in $O(card([X_{q}]_{2\eta}))$. The construction of $S_{p}\parallel S_{q}$ is executed in $O(card([X_{p}]_{2\eta})\cdot card([U_{p}]_{2\mu}) \cdot card([X_{q}]_{2\eta}))$. The construction of $Nb(S_{p}\parallel S_{q})$ is executed in $O(card([X_{p}]_{2\eta}\cap [X_{q}]_{2\eta}) \cdot card([U_{p}]_{2\mu}))$. Hence, the result follows.
%Performing each step of the procedure in Algorithm \ref{alg3} can result in general, in being rather demanding from the complexity point of view. When referring to complexity we consider \textit{space complexity}, related to the size of the implementation of the controller\footnote{Here, the size of $Nb(C^{\ast})$ is interpreted as number of states of $Nb(C^{\ast})$.} $Nb(C^{\ast})$ and
%%\item
%\textit{time complexity}, related to the number of basic algorithmic steps needed to implement $Nb(C^{\ast})$.
%%\end{itemize}
%In this paper we focus on space--complexity issues. The main motivation in this choice is that if the space complexity of algorithms is too much large, since the memory resourses of any computer machine are limited by its nature, the construction of the controller $Nb(C^{\ast})$ could end up in an ``out of memory'' which in turn, would translate in not having a solution to Problem \ref{problem}.
The procedure illustrated in Algorithm \ref{alg3} is not efficient from the space and time complexity point of view\footnote{This qualitative claim will be substantiated in terms of complexity analysis in the next section.} because:
\begin{itemize}
\item It considers the whole state spaces of the plant $P$ and the specification $Q$. A more efficient algorithm would consider only the intersection of the accessible parts of $P$ and $Q$.
\item For any source state $x$ and target state $y$ it includes all transitions $(x,u,y)$ with any control input $u$ by which state $x$ reaches state $y$. A more efficient algorithm would consider for any source state $x$ and target state $y$ only one control input $u$ and hence, only one transition.
\item It first construct the symbolic models $S_{p}$ and $S_{q}$, then the composed system $S_{p}\parallel_{0}S_{q}$ to finally eliminate blocking states from $S_{p}\parallel_{0}S_{q}$. A more efficient algorithm would eliminate blocking states as soon as they show up.
\end{itemize}

Inspired from the research line in the context of on--the--fly verification and control of timed or untimed transition systems (see e.g. \cite{onthefly3,onthefly2}), we now present an algorithm which \textit{integrates each step of the four sub--algorithms in Algorithm \ref{alg3} in only one algorithm}. \\
The proposed procedure is composed of Algorithm \ref{alg} and Algorithm \ref{alg2}. Algorithm \ref{alg} is the main one while Algorithm \ref{alg2} introduces Function \textbf{NonBlock}, which is recursively used in Algorithm \ref{alg}. The outcome of Algorithm \ref{alg} is the symbolic controller $C^{\ast\ast}$ which will be shown in the further results to solve Problem \ref{problem}. 
%We stress that the algorithms that we now propose are designed under the implicit assumption that the plant control system $P$ and the specification system $Q$, together with their symbolic models $S_{p}$ and $S_{q}$ are deterministic. The plant $P$ and the specification $Q$ are deterministic, under the classical assumptions chosen on the vector fields $f_{p}$ and $f_{q}$. Furthermore, the symbolic models $S_{p}$ and $S_{q}$ are shown in Proposition \ref{prop:determinism} to be  deterministic. \\
Given a set $T\subseteq X \times U \times Y$, the set $\mathbf{X}_{source}(T)\subseteq X$ %and $\mathbf{X}_{target}(T)\subseteq Y$%
denotes the projection of $T$ onto $X$, i.e.
\[
\mathbf{X}_{source}(T)=\{x\in X:\exists y\in Y \wedge \exists u\in U \text{ s.t. } (x,u,y)\in T \}.
\]
Given a vector $x\in\mathbb{R}^{n}$ and a precision $\eta\in\mathbb{R}^{+}$, the symbol $[x]_{2\eta}$ denotes the unique vector in $[\mathbb{R}^{n}]_{2\eta}$ such that $x\in\mathcal{B}_{[\eta[}([x]_{2\eta})$.
Algorithm \ref{alg} proceeds as follows. The set of states $X_{0}$ of $C^{\ast\ast}$ is initialized to be $[X_{p,0}\cap X_{q,0}]_{2\eta}$ in line 2.8 and the set of states to be processed, denoted by $\mathbf{X}_{target}$, is initialized to the set of initial states in line 2.9. The set $T$ of transitions and the set $Bad$ of blocking states of $C^{\ast\ast}$ are initialized to be the empty--sets (lines 2.10, 2.11). At each basic step, Algorithm \ref{alg} processes a (non processed) state in line 2.13, by computing the state $y=[\xi^{q}_{x}(\tau)]_{2\eta}$ (line 2.14). If the state $y$ is non--blocking (line 2.15), the algorithm looks for a control input $u\in [U]_{2\mu}$ such that the plant $P$ meets the specification $Q$, i.e. $z=y$ (line 2.20). If such a control input $u$ exists, then boolean variable $Flag$ is updated to $1$ (line 2.21), the transition $(x,u,y)$ is added to the set of transitions $T$ (line 2.25), and the state $y$ is added to the set of the to--be--processed states (line 2.26). If either state $y$ is blocking or no inputs are found for the plant $P$ to meet the specification $Q$, then state $x$ is declared blocking, and Function \textbf{NonBlock}$(T,x,Bad)$ in Algorithm \ref{alg2} is invoked (line 2.30), in order to remove all blocking states originating from $x$. Algorithm \ref{alg} proceeds with further basic steps, until there are no more states to be processed. When Algorithm \ref{alg} terminates, it returns in line 2.34 the symbolic controller $C^{\ast\ast}$. %:
%\[
%C^{\ast\ast}=(\mathbf{X}_{source}(T),X_{0}\cap \mathbf{X}_{source}(T),[U_{p}]_{2\mu},T,Y_{\tau,\eta,\mu},H_{\tau,\eta,\mu}).
%\label{C**}
%\]
Function \textbf{NonBlock}$(T,x,Bad)$ extracts the non--blocking part of $T$. %This algorithm is quite standard in computer science. However, we briefly describe it hereafter for the sake of completeness.
The set $Badx$ includes the states to be processed and is initialized to contain the only state $x$ (line 3.3). At each basic step, for any $y\in Badx$, Function \textbf{NonBlock} removes from the set $T$ any transition $(z,u,y)$ ending up in $y$ (line 3.7), it adds $z$ to the set $Badx$ of states to be processed (line 3.8) and adds $y$ to the set $Bad$ of blocking states (lines 3.11, 3.12). Function \textbf{NonBlock} terminates when there are no more states to be processed and returns in line 2.14 the updated sets of transitions $T$ of and blocking states $Bad$.
Termination of Algorithm \ref{alg} is discussed in the following result:

\incmargin{1em}
\restylealgo{boxed}\linesnumbered
\begin{algorithm}%[H]
\label{alg}
\SetLine
\caption{Integrated Symbolic Control Design.}
\textbf{Input:}\\
%\text{time--delay system $\Sigma=(X,\mathcal{X},\xi_{0},U,\mathcal{U},f)$ satisfying assumptions (A.1-5)}\;
Plant: $P=(X_{p},X_{p,0},U_{p}, \mathcal{U}_{p},f_{p})$\;
Specification: $Q=(X_{q},X_{q,0},U_{q}, \mathcal{U}_{q},f_{q})$\;
Precision: $\varepsilon\in\mathbb{R}^{+}$\;
Parameters: $\theta_{p},\theta_{q}\in\mathbb{R}^{+}$ satisfying (\ref{cond1111})\;
Parameters: $\tau,\eta,\mu\in \mathbb{R}^{+}$ satisfying (\ref{condAUT2})\;
\textbf{Init:}\\
$X_{0}:=[X_{p,0}\cap X_{q,0}]_{2\eta}$\;
$\mathbf{X}_{target}=X_{0}$\;
$T:=\varnothing$\;
$Bad:=\varnothing$\;
%$Term:=0$\;
%\While{$Term=0$}
%{
%$TT=\varnothing \times \varnothing \times \varnothing$\;
%$Term:=1$\;

\ForEach{$x \in [X_{p} \cap X_{q}]_{2\eta}$}
{
\If{$x \in \mathbf{X}_{target} \backslash (\mathbf{X}_{source}(T) \cup Bad)$}
{
           %\If{$x\notin \mathbf{X}_{source}(T)\cup Bad $}
           %{
           %$Term:= 0$\;
           \textbf{compute} $y=[\xi^{q}_{x}(\tau)]_{2\eta}$\;
           \If{$y\notin Bad $}
           %\if{$y\notin Bad$}
           {
           %\eIf{$y\in \mathbf{X}(T)$}
%{
%           $Term:=Term \wedge 1$\;
%           }{
           $Flag:=0$\;
           %     \ForEach{$u \in [U_{p}]_{2\mu}$}
           %     {
           %     \textbf{compute} $z=[\xi^{p}_{xu}(\tau)]_{2\eta}$\;
           %          \If{$z=y$}
           %          {
           %          $Flag:=1$\;
           %          \textbf{break foreach $u \in [U_{p}]_{2\mu}$}\;
           %          }
           %     }                
                \While{$Flag=0$}
                {
                \textbf{choose} $u \in [U_{p}]_{2\mu}$\;
                \textbf{compute} $z=[\xi^{p}_{xu}(\tau)]_{2\eta}$\;
                     \If{$z=y$}
                     {
                     $Flag:=1$\;
                     }                         
                }      
            \If{$Flag=1$}
            {
            $T:=T\cup \{(x,u,y)\}$\;
            $\mathbf{X}_{target}:=\mathbf{X}_{target} \cup \{y\}$\;
            }{
            %$Bad:=Bad\cup\{x\}$\;
            %$(T,Bad):=$\textbf{NonBlock}$(T,x,Bad)$\;
            }
            }{
            \If{$Flag=0 \lor y \in Bad$}
            {$(T,Bad):=$\textbf{NonBlock}$(T,x,Bad)$\;}
            }
            %$\mathbf{X}_{target}=\mathbf{X}_{target} \backslash \{x\}$\;
            }
            }
\textbf{output:} $C^{\ast\ast}=(\mathbf{X}_{source}(T),X_{0}\cap \mathbf{X}_{source}(T),[U_{p}]_{2\mu},T,Y_{\tau,\eta,\mu},H_{\tau,\eta,\mu})$%\;
\end{algorithm}
\decmargin{1em}

\incmargin{1em}
\restylealgo{boxed}\linesnumbered
\begin{algorithm}%[H]
\label{alg2}
\SetLine
\caption{Non--blocking Algorithm.}
\textbf{Function }$(T,Bad):=$\textbf{NonBlock}$(T,x,Bad)$\;
\textbf{Init:}\\
$Badx:=\{x\}$\;
%$Term:=0$\;
%\While{$Term=0$}
%{
%$Term:=1$\;
\ForEach{$y \in Badx$}
{
\ForEach{$z \in \mathbf{X}_{source}(T)$}
{
           \If{$\exists u\in[U]_{2\mu}$ such that $(z,u,y)\in T$}
           {
     			 $T:=T\backslash \{(z,u,y)\}$\;
     			 %\If{$\nexists u\in[U]_{2\mu}$ such that $(z,u,w)\in T$ for some $w$}
           %{     			
     			 $Badx:=Badx\cup \{z\}$\;
     			 %$Term:= 0$\;
     			 %}
           }
}
$Badx:=Badx\backslash \{y\}$\;
$Bad:=Bad\cup \{y\}$\;
}
%}
%$Bad:=Badx\cup Bad$\;
\textbf{output:} $(T,Bad)$
\end{algorithm}
\decmargin{1em}

\begin{theorem}
Algorithm \ref{alg} terminates in a finite number of steps.
\end{theorem}

\begin{proof}
Algorithm \ref{alg} terminates when there are no more states $x$ in $\mathbf{X}_{target}$ to be processed. For each state $x$, either line 2.25 or line 2.30 is executed (depending on the value of the boolean variable $Flag$); this ensures by line 2.13 that state $x$ cannot be processed again in future iterations. Furthermore, the set $\mathbf{X}_{target}$ is nondecreasing (see line 2.26) and always contained in the finite set $[X_{p}]_{2\eta}\cap [X_{q}]_{2\eta}$. Hence, provided that Algorithm \ref{alg2} terminates in finite time, the result follows. Regarding termination of Algorithm \ref{alg2}, in the worst case the set $Bad$ ends up in coinciding with the accessible states of $S_{p}$ and $S_{q}$ (line 3.12) and the set $Badx$ ends up in being empty (line 3.11). Hence from line 3.4, finite termination of Algorithm \ref{alg2} is guaranteed.
\end{proof}

Formal correctness of Algorithm \ref{alg} is guaranteed by the following result.

\begin{theorem}
\label{th1}
Controllers $Nb(C^{\ast})$ and $C^{\ast\ast}$ are exactly bisimilar.
\end{theorem}

\begin{proof}
(Sketch.) For any state $(x_{p},x_{q})$ of the accessible part $Ac(Nb(C^{\ast}))$ of $Nb(C^{\ast})$ there exists a state $x_{c}$ of $C^{\ast\ast}$ so that $x_{p}=x_{q}=x_{c}$ (see lines 2.14, 2.19, 2.20 and 2.25 in Algorithm \ref{alg}). Consider the relation defined by $((x_{p},x_{q}),x_{c})\in \mathcal{R}$ if and only if $x_{p}=x_{c}$. It is readily seen that $\mathcal{R}$ is a $0$--bisimulation relation between $Nb(C^{\ast})$ and $C^{\ast\ast}$.
\end{proof}

By the above result the controller $Nb(C^{\ast})$ solves Problem \ref{problem} if and only if the controller $C^{\ast\ast}$ solves Problem \ref{problem}. Hence, it shows that Algorithm \ref{alg} is correct. While the controllers $Nb(C^{\ast})$ and $C^{\ast\ast}$ are exactly bisimilar, the number of states of $C^{\ast\ast}$ is in general, smaller than the one of $Nb(C^{\ast})$. In fact the controller $Nb(C^{\ast})$ may contain spurious states, e.g. states which are not accessible from a quantized initial condition in $S_{p}$ and a quantized initial condition in $S_{q}$, since in general $Ac(Nb(C^{\ast}))$ is a (strict) sub--system of $Nb(C^{\ast})$. On the other hand, a straightforward inspection of Algorithm \ref{alg} reveals that:

\begin{proposition}
\label{prop}
$Ac(C^{\ast\ast})=C^{\ast\ast}$.
\end{proposition}

Hence, the aforementioned spurious states of $Nb(C^{\ast})$ are not included in $C^{\ast\ast}$. The above remarks suggest the following result:

\begin{theorem}
\label{thmin}
$C^{\ast\ast}$ is the minimal $0$--bisimilar system of $Nb(C^{\ast})$.%, i.e. for any system $C$ which is $0$--bisimilar to $C^{\ast\ast}$, the number of states of $C$ is greather than or equal to the one of $C^{\ast\ast}$.
\end{theorem}

\begin{proof}
The proof can be given by using standard arguments on bisimulation theory \cite{ModelChecking}. Briefly, since by Proposition \ref{prop} $Ac(C^{\ast\ast})=C^{\ast\ast}$ and since the output function $H_{\tau,\eta,\mu}$ of $C^{\ast\ast}$ is the natural inclusion from $\mathbf{X}_{source}(T)$ to $X$, the maximal $0$--bisimulation relation $\mathcal{R}^{\ast}$ between $C^{\ast\ast}$ and itself is the identity relation, i.e. $\mathcal{R}^{\ast}=\{(x_{1},x_{2})\in \mathbf{X}_{source}(T) \times \mathbf{X}_{source}(T) : x_{1}=x_{2}\}$.
%We recall that $\mathcal{R}^{\ast}$ is said to be the maximal $0$--bisimulation relation between $C^{\ast\ast}$ and itself, if for any other $0$--bisimulation relation $\mathcal{R}$ between $C^{\ast\ast}$ and itself, $\mathcal{R}\subseteq \mathcal{R}^{\ast}$. By standard results the minimal bisimilar system $C$ of $C^{\ast\ast}$ is obtained by making the quotient of $C^{\ast\ast}$ induced by the equivalence relation $\mathcal{R}^{\ast}$.
Since $\mathcal{R}^{\ast}$ is the identity relation, the quotient of $C^{\ast\ast}$ induced by $\mathcal{R}^{\ast}$, coincides with $C^{\ast\ast}$. Finally, since by Theorem \ref{th1} systems $C^{\ast\ast}$ and $Nb(C^{\ast})$ are $0$--bisimilar, the result follows.
\end{proof}

The above result is important because it shows that the controller $C^{\ast\ast}$ is the system with the smallest number of states which is equivalent by bisimulation to the solution $Nb(C^{\ast})$ of Problem \ref{problem}.

\section{Space and Time Complexity Analysis}

In this section we provide a formal comparison in terms of space and time complexity analysis, between the procedure illustrated in Algorithm \ref{alg3} and Algorithm \ref{alg}. %We evaluate space complexity, as the number of transitions needed in the construction of the symbolic controllers $Nb(C^{\ast})$ and $C^{\ast\ast}$.

\begin{proposition}
Space complexity of Algorithm \ref{alg3} is $O(\max\{card([X_{p}]_{2\eta}) \cdot card([U_{p}]_{2\mu}),card([X_{q}]_{2\eta})\})$.
\label{SpaceComplexNonIntegrated}
\end{proposition}

\begin{proof}
Since by Proposition \ref{prop:determinism} system $S_{p}$ is deterministic, the number of transitions of $S_{p}$ amounts to $card([X_{p}]_{2\eta}) \cdot card([U_{p}]_{2\mu})$. For the same reason, the number of transitions of $S_{q}$ is given by $card([X_{q}]_{2\eta})$. By definition of exact composition (see Definition \ref{composition} with $\varepsilon=0$), the number of transitions in $S_{p}\parallel_{0} S_{q}$ amounts in the worst case to $(card([X_{p}]_{2\eta}) \cap card([X_{q}]_{2\eta}))\cdot card([U_{p}]_{2\mu})$. By definition of the $Nb$ operator, the number of transitions in $Nb(C^{\ast})$ is less than or equal to the one of $S_{p}\parallel_{0} S_{q}$. Hence, by comparing the above worst case bounds, the result follows.
\end{proof}

\begin{proposition}
Space complexity of Algorithm \ref{alg} is $O(card([X_{p}]_{2\eta}\cap [X_{q}]_{2\eta}))$.
\label{SpaceComplexIntegrated}
\end{proposition}

\begin{proof}
By lines 2.14, 2.15, 2.20, and 2.25 in Algorithm \ref{alg}, the triple $(x,u,y)$ is added to the set $T$ of transitions of $C^{\ast\ast}$, if $(x,u,y)$ is a transition of $S_{p}$ and $(x,y)$ is a transition of $S_{q}$. Hence, the result follows from determinism of systems $S_{p}$ and $S_{q}$, which is guaranteed by Proposition \ref{prop:determinism}.
\end{proof}

By comparing Propositions \ref{SpaceComplexNonIntegrated} and \ref{SpaceComplexIntegrated}, it is readily seen that space complexity of Algorithm \ref{alg} is smaller than or equal to space complexity of Algorithm \ref{alg3}. In particular, when the plant system $P$ and the specification system $Q$ coincide, implying $[X_{p}]_{2\eta}=[X_{q}]_{2\eta}$ and $card([U_{p}]_{2\mu})=1$, the space complexity of the procedure in Algorithm \ref{alg3} and of Algorithm \ref{alg} coincides, resulting in $O(card([X_{p}]_{2\eta}))=O(card([X_{q}]_{2\eta}))$. This is indeed consistent with the integration philosophy that we advocated in Algorithm \ref{alg}. Algorithm \ref{alg} becomes more and more efficient from the space complexity point of view as much as the behaviours of the plant and of the specification differ. When $P$ and $Q$ coincide there is no gain in terms of space complexity, in the use of Algorithm \ref{alg}. We now proceed with a further step by providing a comparison in terms of time complexity analysis.

\begin{proposition}
Time complexity of Algorithm \ref{alg3} is $O(card([ X_{q}]_{2\eta}) \cdot card([X_{p}]_{2\eta}) \cdot card([U_{p}]_{2\mu}))$.
\label{TimeComplexNonIntegrated}
\end{proposition}

\begin{proof}
The number of steps needed in the construction of $S_{p}$ and $S_{q}$ amounts to $card([X_{p}]_{2\eta})\cdot card([U_{p}]_{2\mu})$ and $card([X_{q}]_{2\eta})$, respectively. Since as shown in Proposition \ref{SpaceComplexNonIntegrated} the number of transitions in $S_{p}$ and $S_{q}$ is given respectively by $card([X_{p}]_{2\eta})\cdot card([U_{p}]_{2\mu})$ and $card([X_{q}]_{2\eta})$, the number of steps needed in the
construction of $S_{p}\parallel_{0}S_{q}$ is given by $card([X_{q}]_{2\eta})\cdot card([X_{p}]_{2\eta}) \cdot card([U_{p}]_{2\mu})$. Regarding the computation of the non--blocking part $Nb(S_{p}\parallel_{0} S_{q})$, in the worst case for any state of $S_{p}\parallel_{0} S_{q}$, i.e. for any state in $[X_{q}\cap X_{p}]_{2\eta}$, all transitions in $S_{p}\parallel_{0} S_{q}$ are needed to be processed in order to find blocking states. Since the number of transitions in $S_{p}\parallel_{0} S_{q}$ is $card([X_{q}\cap X_{p}]_{2\eta})\cdot card([U_{p}]_{2\mu})$, the overall number of steps needed in the computation of $Nb(S_{p}\parallel_{0} S_{q})$ is given by $card([X_{q}\cap X_{p}]_{2\eta})^{2}\cdot card([U_{p}]_{2\mu})$. By comparing the above worst case bounds, the result follows.
\end{proof}

\begin{proposition}
Time complexity of Algorithm \ref{alg} is
\[
O(\max \{card([X_{q}\cap X_{p}]_{2\eta}) \cdot card([U_{p}]_{2\mu}),card([X_{q}\cap X_{p}]_{2\eta})^{2}\}).
\]
\label{TimeComplexIntegrated}
\end{proposition}

\begin{proof}
By exploring Algorithm \ref{alg}, it is easy to see that the number of steps needed in the computation of $C^{\ast\ast}$ is upper bounded by:
\begin{equation}
\sum_{i=0}^{N_{1}}(N_{2}+N_{3}),
\label{TC1}
\end{equation}
where $N_{1}=card([X_{p}\cap X_{q}]_{2\eta})$, $N_{2}$ is an upper bound to the number of steps needed in the execution of lines 2.13/27 in Algorithm \ref{alg}, and $N_{3}$ is an upper bound to the number of steps needed in the execution of lines 2.28/30 in Algorithm \ref{alg}. Quantity in (\ref{TC1}) can be rewritten as the sum of the term $\sum_{i=0}^{N_{1}}N_{2}$ and the term $\sum_{i=0}^{N_{1}}N_{3}$, the first of which is upper bounded by $card([X_{p}\cap X_{q}]_{2\eta})\cdot card([U_{p}]_{2\mu})$. Regarding the term $\sum_{i=0}^{N_{1}}N_{3}$, whenever Algorithm \ref{alg} executes line 2.30, i.e. $(T,Bad):=$\textbf{NonBlock}$(T,x,Bad)$, states $x$ involved are different. Indeed suppose by contradiction that at step $i$ state $x$ is processed in line 2.30 and at step $j$ state $x'$ is processed in line 2.30 with $i<j$ and $x=x'$. When at step $i$ Algorithm \ref{alg2} is invoked, state $x$ is added to the set $Bad$ (see lines 3.3, 3.4 and 3.12). Since at the end of step $i$ state $x\in Bad$, in the further steps and in particular at step $j$, state $x$ will be no longer processed (see line 2.13). Since $x'=x$, then at step $j$ state $x'$ cannot be processed in line 2.30. Hence a contradiction holds. Since any time Algorithm \ref{alg2} is invoked it processes different states, the overall time complexity due to the term $\sum_{i=0}^{N_{1}}N_{3}$ is upper bounded by the time complexity needed in computing the non--blocking part of $S_{p}\parallel_{0} S_{q}$ which, from Proposition \ref{TimeComplexNonIntegrated} amounts to $card([X_{q}\cap X_{p}]_{2\eta})^{2}\cdot card([U_{p}]_{2\mu})$. By comparing the above worst case bounds, the result follows.
\end{proof}

By comparing Propositions \ref{TimeComplexNonIntegrated} and \ref{TimeComplexIntegrated}, it is readily seen that time complexity of Algorithm \ref{alg} is smaller than or equal to time complexity of Algorithm \ref{alg3}. In particular, when the plant system $P$ and the specification system $Q$ coincide, implying $[X_{p}]_{2\eta}=[X_{q}]_{2\eta}$ and $card([U_{p}]_{2\mu})=1$, the time complexity of the procedure in Algorithm \ref{alg3} and of Algorithm \ref{alg} coincides, resulting in $O(card([X_{p}]_{2\eta})^{2})=O(card([X_{q}]_{2\eta})^{2})$.

\section{Examples}\label{exampleXXX}
In this section we present some examples of application of the results illustrated in the previous sections. In particular, we consider in Section \ref{Sec:example1} symbolic control design problem for a nonlinear control system 
% with infinite and finite specifications 
and in Section \ref{Sec:example2} symbolic control design for linear control systems. The results shown hereafter are based on computations performed on an Intel Core 2 Duo T5500 1.66GHz laptop with 4 GB RAM.

\subsection{Nonlinear Control Systems}\label{Sec:example1}
%We now present an example of application of the results illustrated in the previous sections. Furthermore, we show experimental results of space and time complexity in the computation of the controllers $Nb(C^{\ast})$ and $C^{\ast\ast}$.  %In particular, we consider in Section \ref{Sec:example1} a symbolic control design problem for a pendulum.

%\label{exampleXXX} 
Consider the following plant nonlinear control system:
\[
P:\left\{
\begin{array}
[l]{l}%
\dot{x}_{1}=-2x_{1}+x_{3}^{2}-u\\
\dot{x}_{2}=2x_{1}-7e^{x_{2}}+7\\
\dot{x}_{3}=-3x_{3}+\frac{3}{4}u^{2},
\end{array}
\right.  \label{example}%
\]
and an infinite states specification, expressed by the following differential equation:
\[
Q:\left\{
\begin{array}
[l]{l}%
\dot{x}_{1}=-3x_{1}+x_{3}^{3}\\
\dot{x}_{2}=x_{1}-5\sin x_{2}\\
\dot{x}_{3}=-x_{2}^{2}-4x_{3}.
\end{array}
\right.
\]
We suppose for simplicity that the plant and the specification systems share the same state space, chosen as: 
\[
X_{p}=X_{q}=[-1,1[\times\lbrack-1,1[\times \lbrack-1,1[,
\]
the same set of initial states, chosen as: 
\[
X_{p}^{0}=X_{q}^{0}=[-1,0[\times\lbrack-1,0[\times\lbrack-1,0[,
\]
and that the plant input space is: 
\[
U=[-1,1].
\]
By using the $\delta$--ISS Lyapunov characterization in \cite{IncrementalS} it is possible to show the plant system $P$ is $\delta$--ISS with functions: 
\[
\begin{array}
{ccc}
\beta_{p}(r,s):=\sqrt{2}\,e^{-1.21\,s}\,r, &
\gamma(r):=\sqrt{14.88\,r}, &
r,s\in\mathbb{R}_{0}^{+}.
\end{array}
\]
Analogously the specification system $Q$ can be shown to be $\delta$--ISS with function: 
\[
\begin{array}
{cc}
\beta_{q}(r,s):=\sqrt{2}\,e^{-s}\,r, &
r,s\in\mathbb{R}_{0}^{+}.
\end{array}
\]
For a precision $\varepsilon=0.2$, we can choose the following quantization parameters for the plant and the specification systems:
\[
\begin{array}{ccccc}
\theta_{p}=0.13, & \theta_{q}=0.07, & \eta=1/30, & \tau=1, & \mu=0.001.
\end{array}
\]
The above choice of quantization parameters guarantees that the inequalities in (\ref{condAUT2}) and (\ref{cond1111}) are fulfilled. By running Algorithm \ref{alg} the integrated symbolic controller $C^{\ast\ast}$ has been designed. Given the large size of the controller obtained ($3152$ states) we do not report in the paper further details on it. Figures \ref{fig1} shows the evolution of the plant system $P$ when interconnected with the symbolic controller $C^{\ast\ast}$ and the evolution of the specification system $Q$, with initial condition $x_{0}=(-1,-1,-1+4\,\eta)$. It is readily seen from the plots that for the initial condition $x_{0}$ the specification is fulfilled, up to the precision $\varepsilon=0.2$ chosen in this example. 

\begin{figure}
\begin{center}
\includegraphics[scale=0.35]{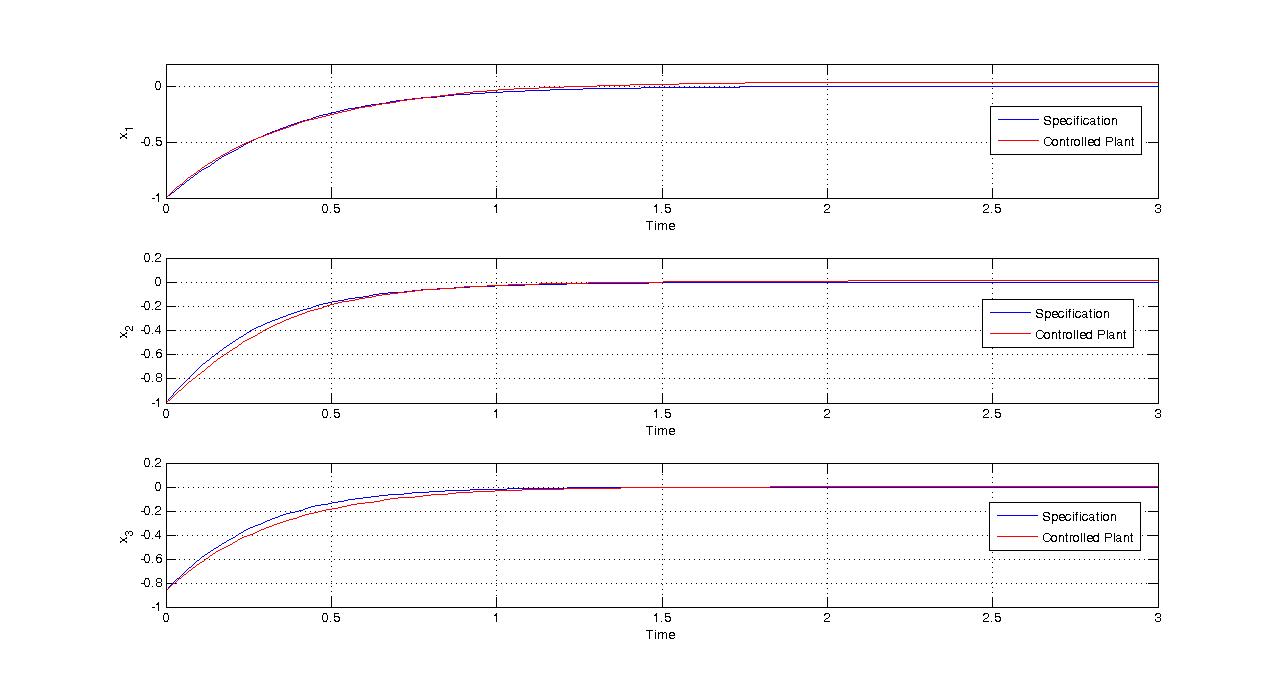}
\caption{Trajectories of the controlled plant and the specification systems with initial condition $(-1,-1,-1+4\,\eta)$.} 
\label{fig1}
\end{center}
\end{figure}

We conclude this section by discussing a comparison between the ``integrated'' approach formulated in Algorithm \ref{alg} and the ``non--integrated'' approach described in Algorithm \ref{alg3}. Experimental results associated with the computation of $C^{\ast\ast}$ and of $Nb(C^{\ast})$ are reported in Tables 1 and 2. In particular, Table 1 shows details in the computation of the controller $Nb(C^{\ast})$ performed by running Algorithm \ref{alg3}. Table 2 reports a comparison between the computation of the controllers $C^{\ast\ast}$ and $Nb(C^{\ast})$. The computation time needed in the construction of the controllers is expressed in seconds and the maximal memory occupation is given in terms of the maximal number of data needed in the construction of the controllers. In particular, the maximal memory occupation in the construction of $Nb(C^{\ast})$ is expressed as the sum of the number of transitions of $S_{p}$, the number of transitions of $S_{q}$ and the number of transitions of $S_{p}\parallel S_{q}$, while the maximal memory occupation in the construction of $C^{\ast\ast}$ is given as the sum of the number of transitions in $C^{\ast\ast}$ and the number of states in $Bad$. For both controllers $Nb(C^{\ast})$ and $C^{\ast\ast}$ each transition is weighted as three data and each state as one datum. The experimental results shown in Table 2 can be summarized, as follows:

\begin{itemize}
\item The number of states of $C^{\ast\ast}$ is $14\%$ times the number of states of $Nb(C^{\ast})$;
%\item The number of blocking states of $C^{\ast\ast}$ is $12\%$ times the number of blocking states of $Nb(C^{\ast})$;
\item The number of transitions of $C^{\ast\ast}$ is $0.25\dot \%$ times the number of transitions of $Nb(C^{\ast})$; 
\item The maximal memory occupation of $C^{\ast\ast}$ is $0.011\%$ times the maximal memory occupation of $Nb(C^{\ast})$; 
\item The time needed in the computation of $C^{\ast\ast}$ is $8\%$ times the time of computation of $Nb(C^{\ast})$. 
\end{itemize}

%We stress that the computation of the system $Nb(C^{\ast})$ is, in general, more expensive because it needs the preliminary
%computation of the symbolic systems $S_{p}$ and $S_{q}$.

\begin{table}[ptb]
%[ptb]%
\begin{tabular}
[c]{lrrrrr}\hline
& $S_{p}$ & $S_{q}$ & $C^{\ast}$ & $Nb(C^{\ast})$ \\
\hline
States & $29791$ & $29791$ & $21894$ & $21894$ \\
Transitions & $29820791$ & $29791$ & $1265217$ & $1265217$ \\
\hline
\end{tabular}
\caption{Details on the computation of $Nb(C^{\ast})$.}%
\label{Tab1}%
\end{table}

\begin{table}[ptb]
%[ptb]%
\begin{tabular}
[c]{lrrrr}\hline
%\textbf{Statistics}
& $Nb\left(  C^{\ast}\right)  $ & $C^{\ast\ast}$ & Ratio & \\\hline
States & $21894$ & $3152$ & $0.14$ & \\
%Blocking states & $7897$ & $944$ & $0.12$ & \\
Transitions & $1265217$ & $3152$ & $2.5\cdot10^{-3}$ & \\
Max memory occupation & $93347397$ & $10400$ & $1.11\cdot10^{-4}$ &
\\
Time & $147487$ & $11144$ & $0.08$ & \\\hline
\end{tabular}
\caption{Comparison in the computation of $Nb(C^{\ast})$ and $C^{\ast\ast}.$}%
\label{Tab2}%
\end{table}

\subsection{Linear Control Systems}\label{Sec:example2}
In this section we consider eight examples randomly chosen in the class of linear systems, characterized by different properties regarding controllability and eigenvalues of dynamical matrices. We consider controllable, versus noncontrollable plant systems (Examples no. 1, 2, 3, 4 vs. 5, 6, 7, 8), plant dynamical matrices $A_{p}$ with real, versus complex eigenvalues (Examples no. 3, 4, 5, 6, 7, 8 vs. 1, 2), specification dynamical matrices $A_{q}$ with real, versus complex eigenvalues (Examples no. 2, 3, 7 vs. 1, 4, 5, 6, 8).
Table 3 shows the experimental results. In particular lines 1.1, 1.2 and 1.3 show respectively, dynamical matrices $A_{p}$ and $B_{p}$ of the plant $P$ and dynamical matrices $A_{q}$ of the specification $Q$. For simplicity we consider in the eight examples the same state space of the plant and the specification, chosen as: 
\[
X_{p}=X_{q}=[-0.5,0.5[\times [-0.5,0.5[,
\]
the same set of initial states of $P$ and $Q$, chosen as: 
\[
X_{p}^{0}=X_{q}^{0}=[-0.25,0.25[\times [-0.25,0.25[, 
\]
and the same input space, chosen as: 
\[
U=[-2,2].
\]
The quantization parameters in the construction of the symbolic systems $S_{p}$ and $S_{q}$ are the same in all the examples and chosen as:
\[
\begin{array}{cccccc}
\varepsilon=0.1, & \tau=0.5, & \mu=0.001, & \eta=0.01, & \theta_{p}=0.05, & \theta_{q}=0.05.
\end{array}
\]
It is readily seen that the above parameters satisfy the inequalities in (\ref{condAUT2}) and (\ref{cond1111}). Experimental results associated with the computation of the controller $Nb(C^{\ast})$ are reported in lines 2.1/2.10. In particular, line 2.10 shows the time of computation needed in the construction of $Nb(C^{\ast})$ and line 2.9 shows the maximal memory occupation in the construction of $Nb(C^{\ast})$. 
Experimental results associated with the computation of the controller $C^{\ast\ast}$ are reported in lines 3.1/3.5. In particular line 3.5 shows the time of computation needed in the construction of $C^{\ast\ast}$ and line 3.4 shows the maximal memory occupation in the construction of $C^{\ast\ast}$. Table 4 summarizes the results shown in Table 3:
\begin{itemize}
\item \textbf{Line 4.1: Gain in terms of number of states.} The minimum gain of the integrated procedure versus the non--integrated procedure is obtained in Example $\#$ $5$, resulting in $100\%$ (meaning that in this example there is no gain in the integrated procedure) and, the maximum gain is obtained in Example $\#$ $7$, resulting in $53\%$.
\item \textbf{Line 4.2: Gain in terms of number of transitions.} The minimum gain of the integrated procedure versus the non--integrated procedure is obtained in Example $\#$ $3$, resulting in $5\%$ and, the maximum gain is obtained in Example $\#$ $7$, resulting in $2\%$.
\item \textbf{Line 4.3: Gain in terms of maximal memory occupation.} The minimum gain of the integrated procedure versus the non--integrated procedure is obtained in Examples $\#$ $2$ and $4$, resulting in $0.017\%$ and, the maximum gain is obtained in Example $\#$ $8$, resulting in $0.007\%$.
\item \textbf{Line 4.4: Gain in terms of time of computation.} The minimum gain of the integrated procedure versus the non--integrated procedure is obtained in Example $\#$ $3$, resulting in $28\%$ and, the maximum gain is obtained in Example $\#$ $7$, resulting in $9\%$.
\end{itemize}

\begin{table}
\begin{tabular}
[c]{lrrrr}
\hline
 & \textbf{Example \# 1} & \textbf{Example \# 2} & \textbf{Example \# 3} & \textbf{Example \# 4} \\
\hline
\textbf{1. Data} & & & & \\
1.1 $A_{p}$ &
$\left(
\begin{array}{cc}
-1 & -0.5 \\
0.5 & -1
\end{array}
\right)$
&
$\left(
\begin{array}{cc}
-0.8 & -0.3 \\
0.3 & -0.8
\end{array}
\right)$
&
$\left(
\begin{array}{cc}
-0.6 & 0 \\
0 & -0.1
\end{array}
\right)$
&
$\left(
\begin{array}{cc}
-0.9 & 0 \\
0 & -0.6
\end{array}
\right)$\\
1.2 $B_{p}$ &
$\left(
\begin{array}{cc}
1 & 1
\end{array}
\right)'$
&
$\left(
\begin{array}{cc}
1 & 1
\end{array}
\right)'$
&
$\left(
\begin{array}{cc}
1 & 1
\end{array}
\right)'$
&
$\left(
\begin{array}{cc}
1 & 1
\end{array}
\right)'$\\
1.3 $A_{q}$ &
$\left(
\begin{array}{cc}
-0.75 & -0.25 \\
0.25 & -0.75
\end{array}
\right)$
&
$\left(
\begin{array}{cc}
-1 & 0 \\
0 & -1
\end{array}
\right)$
&
$\left(
\begin{array}{cc}
-1 & 0 \\
0 & -2
\end{array}
\right)$
&
$\left(
\begin{array}{cc}
-0.8 & -0.4 \\
0.4 & -0.8
\end{array}
\right)$\\
\hline
\textbf{2. }$\mathbf{Nb(C^{\ast})}$  & & & & \\
2.1 States of $S_{p}$  & $2601$ & $2601$ & $2601$ & $2601$ \\
2.2 Transitions of $S_{p}$  & $2675069$ & $2494785$ & $2489327$ & $2446901$ \\
2.3 States of $S_{q}$  & $2601$ & $2601$ & $2601$ & $2601$ \\
2.4 Transitions of $S_{q}$  & $2601$ & $2601$ & $2601$ & $2601$ \\
2.5 States of $C^{\ast}$  & $611$ & $603$ & $403$ & $915$ \\
2.6 Transitions of $C^{\ast}$  & $8013$ & $7507$ & $4969$ & $11151$ \\
2.7 States of $Nb(C^{\ast})$ & $403$ & $521$ & $343$ & $499$ \\
2.8 Transitions of $Nb(C^{\ast})$ & $5719$ & $6753$ & $4331$ & $6505$ \\
2.9 Max($Nb(C^{\ast})$) & $8057049$ & $7514679$ & $7490691$ & $7381959$ \\
2.10 Time($Nb(C^{\ast})$) & $7780$ & $7095$ & $4648$ & $4068$ \\
\hline
\textbf{3. }$\mathbf{C^{\ast\ast}}$  & & & & \\
3.1 States of $C^{\ast\ast}$  & $239$ & $281$ & $199$ & $277$ \\
3.2 Transitions of $C^{\ast\ast}$  & $239$ & $281$ & $199$ & $277$ \\
3.3 States in $Bad$ & $490$ & $448$ & $530$ & $452$ \\
3.4 Max($C^{\ast\ast}$) & $1207$ & $1291$ & $1127$ & $1283$ \\
3.5 Time($C^{\ast\ast}$) & $1300$ & $1800$ & $1300$ & $770$ \\
\hline
\hline
 & \textbf{Example \# 5} & \textbf{Example \# 6} & \textbf{Example \# 7} & \textbf{Example \# 8} \\
\hline
\textbf{1. Data} & & & & \\
1.1 $A_{p}$ &
$\left(
\begin{array}{cc}
-0.9 & 0 \\
0 & -0.6
\end{array}
\right)$
&
$\left(
\begin{array}{cc}
-1.5 & 1 \\
0 & -1.5
\end{array}
\right)$
&
$\left(
\begin{array}{cc}
-0.8 & 0 \\
0 & -0.6
\end{array}
\right)$
&
$\left(
\begin{array}{cc}
-1.5 & 1 \\
0 & -1.5
\end{array}
\right)$\\
1.2 $B_{p}$ &
%$\left(
%\begin{array}{cc}
%x_{11} & x_{12}
%\end{array}
%\right)'$
%&
$\left(
\begin{array}{cc}
1 & 0
\end{array}
\right)'$
&
$\left(
\begin{array}{cc}
1 & 0
\end{array}
\right)'$
&
$\left(
\begin{array}{cc}
1 & 0
\end{array}
\right)'$
&
$\left(
\begin{array}{cc}
1 & 0
\end{array}
\right)'$\\
1.3 $A_{q}$ &
$\left(
\begin{array}{cc}
-0.8 & -0.4 \\
0.4 & 0.8
\end{array}
\right)$
&
$\left(
\begin{array}{cc}
-0.9 & -0.5 \\
0.5 & -0.9
\end{array}
\right)$
&
$\left(
\begin{array}{cc}
-1.2 & 0 \\
0 & -0.5
\end{array}
\right)$
&
$\left(
\begin{array}{cc}
-1 & -0.75 \\
0.75 & -1
\end{array}
\right)$\\
\hline
\textbf{2. }$\mathbf{Nb(C^{\ast})}$  & & & & \\
2.1 States of $S_{p}$  & $2601$ & $2601$ & $2601$ & $2601$ \\
2.2 Transitions of $S_{p}$  & $3290367$ & $3721269$ & $3215397$ & $3721269$ \\
2.3 States of $S_{q}$  & $2601$ & $2601$ & $2601$ & $2601$ \\
2.4 Transitions of $S_{q}$  & $2601$ & $2601$ & $2601$ & $2601$ \\
2.5 States of $C^{\ast}$  & $381$ & $325$ & $1377$ & $227$ \\
2.6 Transitions of $C^{\ast}$  & $9467$ & $9233$ & $33453$ & $6451$ \\
2.7 States of $Nb(C^{\ast})$ & $99$ & $129$ & $153$ & $65$ \\
2.8 Transitions of $Nb(C^{\ast})$ & $2461$ & $3665$ & $3717$ & $1847$ \\
2.9 Max($Nb(C^{\ast})$) & $9907305$ & $11199309$ & $9754353$ & $11190963$ \\
2.10 Time($Nb(C^{\ast})$) & $6285$ & $7080$ & $4880$ & $9444$ \\
\hline
\textbf{3. }$\mathbf{C^{\ast\ast}}$  & & & & \\
3.1 States of $C^{\ast\ast}$  & $99$ & $109$ & $81$ & $53$ \\
3.2 Transitions of $C^{\ast\ast}$  & $99$ & $109$ & $81$ & $53$ \\
3.3 States in $Bad$ & $630$ & $620$ & $648$ & $676$ \\
3.4 Max($C^{\ast\ast}$) & $927$ & $947$ & $891$ & $835$ \\
3.5 Time($C^{\ast\ast}$) & $920$ & $850$ & $430$ & $990$ \\
\hline
\end{tabular}
\caption{Details on the computation of $Nb(C^{\ast})$ and $C^{\ast\ast}$.}
\label{Tab44}
\end{table}

\begin{table}
\begin{tabular}
[c]{lrrrrrrrr}
%\hline
%%\textbf{Statistics} 
%& & & &  \\
%\hline
\hline
 & \textbf{Example \# 1} & \textbf{Example \# 2} & \textbf{Example \# 3} & \textbf{Example \# 4}  \\
\hline
4.1 States($C^{\ast\ast}$)/States($Nb(C^{\ast})$) & $0.59$ & $0.54$ & $0.58$ & $0.55$ \\
4.2 Transitions($C^{\ast\ast}$)/Transitions($Nb(C^{\ast})$) & $0.04$ & $0.04$ & $0.05$ & $0.04$  \\
4.3 Max($C^{\ast\ast}$)/Max($Nb(C^{\ast})$) & $1.5 \cdot 10^{-4}$ & $1.7 \cdot 10^{-4}$ & $1.5 \cdot 10^{-4}$ & $1.7 \cdot 10^{-4}$ \\
4.4 Time($C^{\ast\ast}$)/Time($Nb(C^{\ast})$) & $0.17$ & $0.25$ & $0.28$ & $0.19$ \\
\hline
\hline
 & \textbf{Example \# 5} & \textbf{Example \# 6} & \textbf{Example \# 7} & \textbf{Example \# 8} \\
\hline
4.1 States($C^{\ast\ast}$)/States($Nb(C^{\ast})$) & $1.00$ & $0.84$ & $0.53$ & $0.81$ \\
4.2 Transitions($C^{\ast\ast}$)/Transitions($Nb(C^{\ast})$) & $0.04$ & $0.03$ & $0.02$ & $0.03$ \\
4.3 Max($C^{\ast\ast}$)/Max($Nb(C^{\ast})$) & $0.9 \cdot 10^{-4}$ & $0.8 \cdot 10^{-4}$ & $0.9 \cdot 10^{-4}$ & $0.7 \cdot 10^{-4}$ \\
4.4 Time($C^{\ast\ast}$)/Time($Nb(C^{\ast})$) & $0.15$ & $0.12$ & $0.09$ & $0.10$ \\
\hline
\end{tabular}
\caption{Comparison between the computation of $Nb(C^{\ast})$ and of $C^{\ast\ast}$.}
\label{Tab33}
\end{table}

\section{Discussion}
In this paper we addressed the problem of symbolic control design of nonlinear systems with infinite states specifications, modelled by differential equations. %The solution $NB(C^{\ast})$ of Problem \ref{problem} has been shown to be the non--blocking part of the exact parallel composition of the symbolic systems associated with the plant and the specification. 
After having provided an explicit solution to the symbolic control design problem, we presented Algorithm \ref{alg} which integrates the design of the symbolic controller with the construction of the symbolic systems of the plant and of the specification. 
%????? The results developed in this paper follow the research line concerning symbolic control design of continuous and hybrid systems, as investigated in \cite{LTLControl,BeltaCDC09,TabuadaTAC08} and in particular in the direction of \cite{TabuadaTAC08}. A comparison with the work in \cite{TabuadaTAC08} is discussed hereafter:
%\begin{itemize}
%\item While the present work faces the problem of symbolic control design with infinite states specifications, the work in \cite{TabuadaTAC08} considers finite states specifications, described by automata.
%\item While the present work considers from the early beginning of the control design process non--blocking behaviours, the work in \cite{TabuadaTAC08} do not.
%\item While the present work proposes an integrate scheme for symbolic control design, the work in \cite{TabuadaTAC08} do not.
%\end{itemize}
%????? 
Although the focus of the present paper is on infinite states specifications, it can be shown that the results here presented can be easily adapted to consider finite states specifications which include language specifications, formalized through automata theory \cite{CassandrasBook}. 
This is important because, as shown in the work of \cite{LTLControl,TabuadaTAC08,Belta:06}, automata theory provides a novel class of specifications which were traditionally not addressed before, in the control design of continuous (nonlinear) systems. 
%The symbolic controller $C^{\ast\ast}$, outcome of Algorithm \ref{alg} was shown in Theorem \ref{thmin} to be the minimal exactly bisimilar system of the controller $Nb(C^{\ast})$ which was shown to solve the controller synthesis problem.
%Space and time complexity in the integrated procedure formalized in Algorithm \ref{alg} has been proved in Section 6 to be smaller than or equal to the non--integrated procedure illustrated in Algorithm \ref{alg3}. %The integrated symbolic control design has been based on a hi--level algorithm whose performance in terms of space and time complexity can be certainly improved by appropriately choosing suitable data structures and software layers. This will be the object of future investigation.
Future work will focus on more efficient techniques at the software layer which can further reduce space and time complexity in the implementation of Algorithm \ref{alg}. Useful insights in this direction can be found in the tool \textsf{Pessoa} \cite{pessoa} which employes binary decision diagrams \cite{CUDD} as data structures to encode symbolic systems.
%More efficient techniques to reduce the space complexity in the computation of the controller $C^{\ast\ast}$ can be indeed considered, by designing suitable data structures at the software implementation layer, which could further reduce the memory occupation of computer machines. By taking inspiration from the tool \textsf{Pessoa} \cite{pessoa}, a promising approach could be the use of Binary Decision Diagrams \cite{CUDD}, as appropriate data structures to encode the symbolic controller $C^{\ast\ast}$. This is one of the research issues which will be addressed in further work.

\textbf{Acknowledgement.} The first author would like to thank Paulo Tabuada for having inspired the idea of integration of control algorithms with the construction of the symbolic systems of the plant and the specification.

%\end{acknowledgement}

\section*{Appendix: Notation}\label{appendix}

The identity map on a set $A$ is denoted by $1_{A}$. Given two sets $A$ and $B$, if $A$ is a subset of $B$ we denote by $1_{A}:A\hookrightarrow B$ or
simply by $\imath$ the natural inclusion map taking any $a\in A$ to $\imath\left(  a\right)  =a\in B$. Given a function $f:A\rightarrow B$ the symbol $f(A)$ denotes the image of $A$ through $f$, i.e. $f(A):=\{b\in B:\exists a\in A$ s.t. $b=f(a)\}$; if $C\subset A$ we denote by $f|_{C}$ the restriction of $f$ to $C$, i.e. $f|_{C}(x):=f(x)$ for any $x\in C$. Given a relation $R\subseteq A\times B$, $R^{-1}$ denotes the inverse relation of $R$, i.e. $R^{-1}:=\{\left(  b,a\right)  \in B\times A:\left(  a,b\right)  \in A\times B\}$. A relation $R\subseteq A\times B$ is a preorder if it is reflexive, transitive but not symmetric. The symbols $\mathbb{N}$, $\mathbb{Z}$, $\mathbb{R}$, $\mathbb{R}^{+}$ and $\mathbb{R}_{0}^{+}$ denote the set of natural, integer, real, positive real, and nonnegative real numbers, respectively. Given a vector $x\in\mathbb{R}^{n}$, we denote by $x_{i}$ the $i$--th element of $x$ and by $\Vert x\Vert$ the infinity norm of $x$, we recall that \mbox{$\Vert x\Vert=\max\{|x_1|,|x_2|,...,|x_n|\}$}, where $|x_{i}|$ denotes the absolute value of $x_{i}$. 
%A function \mbox{$f:[a,b]\rightarrow\mathbb{R}^n$} is said to be absolutely continuous on $[a,b]$ if for any $\varepsilon\in\mathbb{R^+}$ there exists $\delta\in\mathbb{R}^+$ so that for every $k\in\mathbb{N}$ and for every sequence of points \mbox{$a\leq{a_1}<b_1<a_2<b_2<\ldots<a_k<b_k\leq{b}$}, if \mbox{$\sum_{i=1}^m(b_i-a_i)<\delta$} then \mbox{$\sum_{i=1}^m\Vert f(b_i)-f(a_i)\Vert<\varepsilon$}. A function \mbox{$f:]a,b[\rightarrow\mathbb{R}^n$} is said to be locally absolutely continuous if the restriction of $f$ to any compact subset of $]a,b[$ is absolutely continuous.
Given a measurable function
\mbox{$f:\mathbb{R}_{0}^{+}\rightarrow\mathbb{R}^n$}, the (essential) supremum
of $f$ is denoted by $\Vert f\Vert_{\infty}$; we recall that $\Vert
f\Vert_{\infty}=(\operatorname*{ess})\sup\left\{  \Vert f\left(  t\right)
\Vert,t\geq0\right\}  $; $f$ is essentially bounded if $\Vert{f}\Vert_{\infty
}<\infty$.
%For a given time $\tau\in\mathbb{R}^+$, define $f_{\tau}$ so that $f_{\tau}(t)=f(t)$, for any $t\in[0,\tau)$, and $f(t)=0$ elsewhere; $f$ is said to be locally essentially bounded if for any $\tau\in\mathbb{R}^+$, $f_{\tau}$ is essentially bounded.
Given $x\in\mathbb{R}^{n}$ and $\varepsilon\in\mathbb{R}^{+}$, the symbol $\mathcal{B}_{\varepsilon}(x)$ denotes the set $\{x\in\mathbb{R}^{n}: \Vert x \Vert\leq \varepsilon\}$ and the symbol $\mathcal{B}_{[\varepsilon[}(x)$ denotes the set $[-\varepsilon+x_{1},x_{1}+\varepsilon[\times [-\varepsilon+x_{2},x_{2}+\varepsilon[\times ... \times [-\varepsilon+x_{n},x_{n}+\varepsilon[$. It is readily seen that if $x\in \mathcal{B}_{\varepsilon}(y)$ and $y\in \mathcal{B}_{[\theta[}(z)$ then
$x\in \mathcal{B}_{[\varepsilon+\theta[}(z)$.
%$\mathcal{B}_{\varepsilon}(x)=\left[  x_{1}-\varepsilon,x_{1}+\varepsilon\right)\times\left[  x_{2}-\varepsilon,x_{2}+\varepsilon\right)  \times...\times\left[  x_{n}-\varepsilon,x_{n}+\varepsilon\right)  $.
For any
\mbox{$A\subseteq\mathbb{R}^{n}$} and \mbox{$\mu\in{\mathbb{R}^+}$}, define
\mbox{$[A]_{\mu}=\{a\in
A\,\,|\,\,a_{i}=k_{i}\mu,k_{i}\in\mathbb{Z},i=1,2,...,n\}$}. The set
$[A]_{\mu}$ will be used as an approximation of the set A with precision $\mu/2
$.
For a given time $\tau
\in\mathbb{R}^{+}$, define $f_{\tau}$ so that $f_{\tau}(t)=f(t)$, for any
$t\in\lbrack0,\tau[$, and $f(t)=0$ elsewhere; $f$ is said to be locally
essentially bounded if for any $\tau\in\mathbb{R}^{+}$, $f_{\tau}$ is
essentially bounded.
%A function $f:\mathbb{R}^{n}\rightarrow\mathbb{R}$ is said to be radially unbounded if $f(x)\rightarrow\infty$ as $\Vert x\Vert\rightarrow\infty$.
A continuous function
\mbox{$\gamma:\mathbb{R}_{0}^{+}\rightarrow\mathbb{R}_{0}^{+}$}, is said to
belong to class $\mathcal{K}$ if it is strictly increasing and
\mbox{$\gamma(0)=0$}; $\gamma$ is said to belong to class $\mathcal{K}%
_{\infty}$ if \mbox{$\gamma\in\mathcal{K}$} and $\gamma(r)\rightarrow\infty$
as $r\rightarrow\infty$. A continuous function
\mbox{$\beta:\mathbb{R}_{0}^{+}\times\mathbb{R}_{0}^{+}\rightarrow\mathbb{R}_{0}^{+}$}
is said to belong to class $\mathcal{KL}$ if, for each fixed $s$, the map
$\beta(r,s)$ belongs to class $\mathcal{K}_{\infty}$ with respect to $r$ and,
for each fixed $r$, the map $\beta(r,s)$ is decreasing with respect to $s$ and
$\beta(r,s)\rightarrow0$ as \mbox{$s\rightarrow\infty$}.

\bibliographystyle{alpha}
\bibliography{biblio1}

\end{document}